\documentclass [12pt,a4paper,reqno]{amsart}

  \input amssymb.sty

\DeclareMathOperator{\Ln}{ln}

\newcommand{\ef}{\end{equation}}
\chardef\bslash=`\\ 

 \newtheorem{thm}{Theorem}[section]
\newtheorem*{thm*}{Theorem}

 \newtheorem*{lem1.2*}{``Lemma 1.2"}
 \newtheorem*{thm1.5*}{``Theorem 1.5"}

\newtheorem{cor}{Corollary}[thm]

\newtheorem{lem}{Lemma}[section]
\newtheorem{corl}{Corollary}[lem]

\newtheorem{assertion}{Assertion}[section]

\theoremstyle{definition}
\newtheorem{defn}{Definition}[section]
\newtheorem{examp}{Example} [section]
\newtheorem*{remark*}{Remarks}
\newtheorem*{MainCon*}{Main conclusion}
\newtheorem*{defn*}{Definition}

 \theoremstyle{remark}
\newtheorem{remark}{Remark}[section]

\newcommand{\thmref}[1]{Theorem~\ref{#1}}
\newcommand{\secref}[1]{Section~\ref{#1}}

\newcommand{\lemref}[1]{Lemma~\ref{#1}}
\newcommand{\corref}[1]{Corollary~\ref{#1}}
\newcommand{\exampref}[1]{Example~\ref{#1}}

\newcommand{\corlref}[1]{Corollary~\ref{#1}}

\newcommand{\A}{\mathcal{A}}

 \renewcommand{\sectionmark}[1]{}
\renewcommand{\Im}{\operatorname{Im}}
\renewcommand{\Re}{\operatorname{Re}}

\newcommand{\ve}{\varepsilon}

\newcommand{\iy}{\infty}

\newcommand{\g}{\gamma}
 
\newcommand{\Dl}{\Delta}

\newcommand{\ov}{\overline}
\newcommand{\vp}{\varphi}

\renewcommand{\Im}{\operatorname{Im}}

\newcommand{\doe}{\overset{\text{def}}{=}}

 \date{}
\begin{document}

\title[Necessary and sufficient conditions for solvability]  {Necessary and sufficient
conditions for solvability\\ of the Hartman-Wintner problem\\ for difference equations}
\author[N.A. Chernyavskaya and L.A. Shuster]{N.A. Chernyavskaya\\
Department of Mathematics and Computer Science\\
Ben-Gurion University of the Negev\\
P.O.B. 653, Beer-Sheva, 84105, Israel\\
\quad\\
L.A. Shuster\\
 Department of Mathematics\\
 Bar-Ilan University, 52900 Ramat Gan, Israel}
\begin{abstract} The equation
\begin{equation}\label{1}
\Dl(r_{n-1}\Dl y_{n-1})=(q_n+\sigma_n)y_n,\quad n\ge 0
\end{equation}
is viewed as a perturbation of the equation
\begin{equation}\label{2}
\Dl(r_{n-1}\Dl z_{n-1})=q_nz_n,\quad n\ge 0
\end{equation}
which does not oscillate at infinity.
The sequences $\{r_n\}_{n=0}^\infty,$\ $\{q_n\}_{n=0}^\infty$ are assumed   real, $r_n>0$
for all $n\ge0,$ the sequences $\{\sigma_n\}_{n=0}^\infty$ may be complex-valued.
We study the Hartman-Wintner problem on asymptotic ``integration" of \eqref{1} for large $n$
in terms of solutions of \eqref{2} and the perturbation $\{\sigma_n\}_{n=0}^\iy$.
\end{abstract}

\maketitle

\section{Introduction}\label{introduction}
\baselineskip 20pt
\setcounter{equation}{0}
\numberwithin{equation}{section}

In the present paper we consider difference equations
\begin{equation}\label{1.1}
\Dl(r_{n-1}\Dl y_{n-1})=(q_n+\sigma_n)y_n,\quad n=0,1,2,\dots
\end{equation}
\begin{equation}\label{1.2}
\Dl(r_{n-1}\Dl z_{n-1})=q_nz_n,\quad n=0,1,2,\dots
\end{equation}
where the sequences $\{r_n\}_{n=0}^\infty,$\ $\{q_n\}_{n=0}^\infty$ are assumed  real,
$r_n>0$ for all $n\ge0,$ and the sequence $\{\sigma_n\}_{n=0}^\infty$ may be complex-valued; $\Dl a_n=a_{n+1}-a_n,$\
$n\ge0,$ for every sequence $\{a_n\}_{n=0}^\iy.$ We assume that equation \eqref{1.2} does not oscillate at infinity.
It is known \cite[Ch.VI]{1} that in this case  it has solutions
  $\{u_n\}_{n=0}^{\infty}$ (the principal, or recessive, solution) and $\{v_n\}_{n=0}^\infty$ (the non-principal, or
dominant, solution),
 and there is an integer
$n_0$ such that  the following relations hold:
 \begin{equation}
\begin{aligned}\label{1.3}
&u_n>0,\quad v_n>0\quad\text{for}\quad n\ge n_0\\
 &r_n(v_{n+1}u_n-u_{n+1}v_n)=1\quad\text{for}\quad n\ge0 \\
 &\lim_{n\to\infty}\frac{u_n}{v_n}=0,\quad\sum_{n=n_0}^\infty\
  \frac{1}{r_nu_nu_{n+1}}=\infty,\quad \sum_{n=n_0}^\infty\ \frac{1}{r_nv_nv_{n+1}}<\infty.
\end{aligned}
 \end{equation}

Our goal is to study the following problem: To find conditions on $\{\sigma_n\}_{n=0}^\infty$
  under which there exists a fundamental system of solutions
(FSS) $\{\tilde u_n,\tilde v_n\}_{n=0}^\iy$ of equation \eqref{1.1} such that
\begin{equation}\label{1.4}
\lim_{n\to\infty}\frac{\tilde u_n}{u_n}=\lim_{n\to\iy}\frac{\tilde v_n}{v_n}=1
\end{equation}
\begin{equation}\label{1.5}
\frac{\tilde u_{n+1}}{\tilde u_n}=\frac{u_{n+1}}{u_n}+o\left(\frac{1}{r_nu_nv_n}\right)\quad\text{for}\quad n\to\iy
\end{equation}
\begin{equation}\label{1.6}
\frac{\tilde v_{n+1}}{\tilde v_n}=\frac{v_{n+1}}{v_n}+o\left(\frac{1}{r_nu_nv_n}\right)\quad\text{for}\quad n\to\iy.
\end{equation}
Here we assume that a FSS $\{u_n,v_n\}_{n=0}^\iy$ of equation \eqref{1.2} with properties \eqref{1.3} is known.

A similar question for differential equations was first studied by P. Hartman and A. Wintner, and
therefore we relate the above problem to their names (and denote it throughout   as problem
\eqref{1.4}--\eqref{1.6}); see \cite{2},\cite{6},\cite{7}, \cite{8},\cite{10} for generalities on problems
\eqref{1.4}--\eqref{1.6}; see \cite{3} for a summary of results on the Hartman-Wintner problem.
Note that a problem close to \eqref{1.4}--\eqref{1.6} was studied in \cite{11}, and the main result of \cite{11} can
be interpeted as a solution of problem \eqref{1.4}--\eqref{1.6}.
Therefore below we present statements of that paper (``projecting" them on to the Hartman-Wintner problem) and
emphasize that the problem considered there is different from problem \eqref{1.4}--\eqref{1.6}.
Here and  throughout the sequel  we say that problem \eqref{1.4}--\eqref{1.6} is solvable if \eqref{1.1} has a FSS
$\{\tilde u_n,\tilde v_n\}_{n=0}^\iy$ satisfying \eqref{1.4}--\eqref{1.6}; in \eqref{1.3}, without loss of generality,
we assume $n_0$ equal to zero.

\begin{assertion}\label{assert1} (\cite{11})\
If the series
\begin{equation}\label{1.7}
J\doe\sum_{n=0}^\infty \sigma_nu_nv_n
\end{equation}
converges (at least conditionally), then the sequence
\begin{equation}\label{1.8}
C_n\doe\frac{v_n}{u_n}\sum_{k=n}^\infty \sigma_ku_k^2,\quad n=0,1,2,\dots
\end{equation}
is well defined, and the following inequalities hold:

\begin{equation}\label{1.9}
|C_n|\le 2  A_n,\quad n=0,1,\dots;\quad   A_n\doe\sup_{m\ge n}|J_m|,\quad
J_m\doe\sum_{k=m}^\infty \ \sigma_ku_kv_k.
\end{equation}
\end{assertion}

\begin{thm}\label{thm1.1}
(\cite{11})\
Suppose that the series $J$ (see \eqref{1.7}) converges (at least, conditionally), and, in addition,
\begin{equation}\label{1.10}
\sum_{n=0}^\infty\frac{A_{n+1}|C_{n+1}|}{r_nu_nv_{n+1}}<\infty
\end{equation}
\begin{equation}\label{1.11}
\lim_{m\to\infty}\sup_{n\ge
m}\frac{1}{A_n}\sum_{k=n}^\infty\frac{A_{n+1}|C_{n+1}|}{r_nu_nv_{n+1}}=\mu<\frac{1}{2}.
\end{equation}
 Then the Hartman-Wintner
problem is solvable.
\end{thm}

\begin{cor}\label{cor1.1}
If the series $J$ (see \eqref{1.7}) absolutely converges, then problem \eqref{1.4}--\eqref{1.6} is solvable.
\end{cor}

We prove Assertion \ref{assert1} in \S5.
\thmref{thm1.1} is proved in \S8.
Note that assumption \eqref{1.11} in its statement turned out to be superfluous.

Let us now formulate our results.
The following theorem contains necessary conditions for solvability of problem \eqref{1.4}--\eqref{1.6}.

\begin{thm}\label{thm1.2}
If the Hartman-Wintner problem   is solvable, then the series $\sigma:$
\begin{equation}\label{1.12}
\sigma\doe \sum_{n=0}^\iy \sigma_nu_n^2
\end{equation}
converges (at least conditionally), the sequence $\{C_n\}_{n=0}^\infty$ (see \eqref{1.8}) is
well defined, and,  in addition,
\begin{equation}\label{1.13}
\lim_{n\to\infty}C_n=0.
\end{equation}
\end{thm}

For the reader's convenience, let us outline the plan of the paper.
Below we follow the general approach to the study of problem \eqref{1.4}--\eqref{1.6}
which was proposed in \cite{3} for its differential analogue.
This means that we restrict the initial problem in order to obtain a criterion for solvability of the narrow problem
 (see Definition 1.1 and
 \thmref{thm1.3} below), and then  find precise a priori requirements to
$\{\sigma_n\}_{n=0}^\iy$ under which the initial problem is equivalent to the narrow one (see \thmref{thm1.4} below).
Finally, we analyze why the conditions for solvability of the Hartman-Wintner problems for differential and difference
equations are not completely analogous whereas the statements of these problems are completely analogous (see
\cite{3}).
We also present examples to all the main results of the paper (see \S9).
We now go over to precise statements.

\begin{defn}\label{def1.1}
Problem \eqref{1.4}--\eqref{1.6} with the additional requirement
\begin{equation}
\label{1.$*$}
\sum_{n=0}^\infty r_nu_nv_{n+1}\left|\frac{\tilde u_{n+1}}{u_{n+1}}-\frac{\tilde
u_n}{u_n}\right|^2<\infty
\tag{1.$*$}
\end{equation}
is called the narrow Hartman-Wintner problem (and is referred to below as problem
 \eqref{1.4}--\eqref{1.$*$}).
\end{defn}

To give a criterion for solvability of problem \eqref{1.4}--\eqref{1.$*$}  and its consequences, we need the
following series:
\begin{equation}\label{1.14}
G\doe \sum_{n=0}^\infty \frac{|C_{n+1}|^2}{r_nu_nv_{n+1}}\qquad \text{(see \eqref{1.8})}
\end{equation}
\begin{equation}\label{1.15}
L\doe\sum_{n=0}^\infty \frac{Re(J_{n+1}\ov C_{n+1})}{r_nu_nv_{n+1}}\qquad \text{(see \eqref{1.9})}
\end{equation}
\begin{equation}\label{1.17}
P\doe\sum_{n=0}^\iy \Re(\sigma_n\overline C_n)u_nv_n.
\end{equation}

 \begin{lem}\label{lem1.1}
Suppose that the series $J$ (see \eqref{1.7}) converges (at least, conditionally).
Then the series $L$ converges  (at least, conditionally) if and only if the series $G$ converges, and the series $P$
converges (at least, conditionally) if and only if the series $G$ (see \eqref{1.14} and $B$ both converge, where
\begin{equation}\label{1.17a}
B\doe \sum_{n=0}^\iy|\sigma_nu_nv_n|^2.
\end{equation}
\end{lem}

The next theorem contains the main result of the paper.

\begin{thm}\label{thm1.3}
The narrow Hartman-Wintner problem is solvable if and only if the series $J$ converge (at least, conditionally; see
\eqref{1.7})  and either  the series
$G$ or $L$   converge ($L$ at least, conditionally).
\end{thm}

\begin{cor}\label{cor1.3.1}
The narrow Hartman-Wintner problem is solvable provided the series $J$  (see\eqref{1.7}) converges (at least, conditionally) and
any of the following conditions $I)$--$IV)$ holds:
\begin{enumerate}
\item[I)]
at least one of the following inequalities holds (see \eqref{1.9}):
\begin{equation}\label{1.16}
\sum_{n=0} ^\iy
\frac{|J_{n+1}C_{n+1}|}{r_nu_nv_{n+1}}<\iy,\quad\sum_{n=0}^\iy\frac{A_{n+1}|C_{n+1}|}{r_nv_{n+1}u_n}<\iy,\quad
\sum_{n=0}^\iy\frac{|J_{n+1}|^2}{r_nu_nv_{n+1}}<\iy.
\end{equation}
\item[II)] the series $P$ converges (at least, conditionally);
\item[III)] the following inequality holds (see \eqref{1.9}):
\begin{equation}\label{1.18}
\sum_{n=0}^\iy|\sigma_n|A_nu_nv_n<\iy
\end{equation}
\item[IV)] the series $J$ (see \eqref{1.7}) absolutely converges.
\end{enumerate}

\end{cor}

\thmref{thm1.3} contains a criterion for solvability of problem  \eqref{1.4}--\eqref{1.$*$} (but not of problem
\eqref{1.4}--\eqref{1.6}).
Therefore one can pose the problem: to find requirements to the perturbation $\{\sigma_n\}_{n=0}^\iy$ under which the
Hartman-Wintner and its narrow analogue are indistinguishable (equivalent).
Clearly, if problem  \eqref{1.4}--\eqref{1.$*$} is solvable, then problem  \eqref{1.4}--\eqref{1.6} is also solvable.
Therefore the above question can be reduced to the following: to find precise a priori requirements to
$\{\sigma_n\}_{n=0}^\iy$  under which any solution of problem  \eqref{1.4}--\eqref{1.6} (i.e., any FSS $\{\tilde
u_n, \tilde v_n\}_{n=0}^\iy$ of equation \eqref{1.1} satisfying relations  \eqref{1.4}--\eqref{1.6}) also satisfies
condition  \eqref{1.$*$}.
In the last case, the Hartman-Wintner problem is completely reduced to the narrow Hartman-Wintner problem, and we say
that the two problems are equivalent.
A criterion for the equivalence of problem \eqref{1.4}--\eqref{1.6} and \eqref{1.4}--\eqref{1.$*$} is given in
\thmref{thm1.4}.

\begin{thm}\label{thm1.4}
Suppose that the Hartman Wintner problem    is solvable and $\{\tilde
u_n, \tilde v_n\}_{n=0}^\iy$ is a FSS of \eqref{1.1} satisfying \eqref{1.4}--\eqref{1.6}.
Then inequality \eqref{1.$*$} holds if and only if $G<\iy$ (see \eqref{1.14}).
\end{thm}

Note that \thmref{thm1.4} often (but not always,  see examples in \S9) allows one to reduce the study of problem
\eqref{1.4}--\eqref{1.6} to the study of problem \eqref{1.4}--\eqref{1.$*$}, i.e., to the application of
\thmref{thm1.3}.
More precisely, the set of equations \eqref{1.1} for which such a reduction is impossible is not larger than the set
of equations \eqref{1.1} for  which the sequence $\{C_n\}_{n=0}^\iy$ is defined, $C_n\to 0$, as $n\to\iy,$ and $G=\iy$ (see
\eqref{1.14}).

Let us now analyze what is common and what is different in the statement on solvability of the Hartman-Wintner
problems for differential and difference equations.
For brevity, we shall do that by conditional comparison of results, as follows.
We assume that the criteria for solvability of the narrow Hartman-Wintner problems for differential and difference
equations are completely analogous.
Under this assumption, we obtain the following assertions (we put them in quotation  marks because they are false):

\begin{lem1.2*}\label{lem1.2}
Suppose that the series $J$ (see \eqref{1.7}) converges (at least, conditionally).
Then the series $L$ and $P$ (see \eqref{1.15}--\eqref{1.17}) converge (at least, conditionally) if and only if the
series $G$ (see \eqref{1.14}) converges.
\end{lem1.2*}

\begin{thm1.5*}\label{thm1.5}
The narrow Hartman-Wintner problem is solvable if and only if condition $\alpha$) and any of conditions $\beta$),
$\gamma$), or $\theta$) hold:
\begin{enumerate}
\item[$\alpha$)] the series $J$ (see \eqref{1.7}) converges (at least, conditionally);
\item[$\beta$)] the series $G$ (see \eqref{1.14} converges;
\item [$\gamma$)] the series $L$ (see \eqref{1.15} converges (at least, conditionally);
\item[$\theta$)] the series $P$ (see \eqref{1.16}) converges (at least, conditionally).
\end{enumerate}
\end{thm1.5*}

Let us now compare ``Theorem~1.5" with \thmref{1.3}.
We see that in the parts $\alpha$)--$\beta$) and $\alpha$)--$\gamma$), ``Theorem~1.5" is true  (and thus completely
analogous to the corresponding theorem  for differential equations).
Here we see the complete analogy between the narrow Hartman-Wintner problems for differential and difference
equations.
In the part $\alpha$)--$\theta$), ``Theorem~1.5" is no longer true, and here we see the difference between these
problems.
According to \corref{cor1.3.1}, in the part $\alpha$)--$\theta$) ``Theorem~1.5" (difference case) is only  a
sufficient (but not a necessary) condition for solvability of problem \eqref{1.4}--\eqref{1.$*$}).
Thus, ``discretization" of the narrow Hartman-Wintner problem shows that ``Theorem~1.5" is no more precise in the
part $\alpha$)--$\theta$).
The reason can be explained by comparing \lemref{lem1.1} with ``Lemma 1.2".
The convergence (at least, conditional) of both the series $J$ and $P$ (see \eqref{1.7}, \eqref{1.17} implies not
only the convergence of the series $G$ (which already gives a criterion for solvability of the narrow Hartman-Wintner
problem) but also the convergence of the series $B$ (see \eqref{1.17a}) which is, in general,  is not obligatory and
can be viewed as some ``additional load" on the parameters of the problem.

\begin{examp}\label{examp1}
Consider the Hartman-Wintner problem for the equations
\begin{equation}\label{1.20}
\Delta(n^\alpha\Dl y_n)=\frac{(-1)^n}{(n+1)^\beta}\ y_{n+1},\quad n\ge 1
\end{equation}
\begin{equation}\label{1.21}
\Delta(n^\alpha\Dl z_n)=0\cdot z_{n+1},\quad n\ge 1
\end{equation}
for $\alpha\in[0,1)$, $\beta\in\left(1-\alpha,\left.\frac{3}{2}\right.  -\alpha\right].$
In this case (see  Section 9, \exampref{examp9.3}), problem \eqref{1.4}--\eqref{1.$*$}) is solvable but the series $B$ (see
\eqref{1.17a}) diverges and, according to \lemref{lem1.1}, the series $P$ (see \eqref{1.17}) also diverges, which
contradicts the assertion $\alpha$)--$\theta$) of ``Theorem~1.5".
\end{examp}

To conclude our analysis, we now have to find a reason which explains the difference between Lemmas \ref{lem1.1} and
``1.2".
In this case such a difference arises because the integral calculus and the finite difference calculus are not always
completely analogous.
More precisely,  in the case of the Hartman-Wintner problem for differential equations the proof of the criterion for
solvability of type $\alpha$)--$\theta$) (see ``Theorem~1.5") relies, after all, upon the following obvious statement
(see \cite{3}): if $x(t)$ is a continuously differentiable function on $[1,\iy)$ such that $x(t)\to0$ as $t\to\iy,$
then the integral
\begin{equation}\label{1.22}
\int_1^\iy x'(t)x(t)dt
\end{equation}
converges.
The difference analogue of this integral (with step 1) is of the form
\begin{equation}\label{1.23}
\sum_{n=1}^\iy (\Dl x_n)x_n,\quad x_n=x(n),\quad n=1,2,\dots;\quad x_n\to0\quad\text{as}\quad n\to \iy.
\end{equation}
(It is this series that appears in the course of the analysis of the condition $\alpha$)--$\theta$) of
\linebreak ``Theorem~1.5", see \S5).

It is easy to see that the series \eqref{1.23} converges if and only if the series
\begin{equation}\label{1.24}
D=\sum_{n=1}^\iy (\Dl x_n)^2
\end{equation}
converges.
Indeed,
\begin{equation}\label{1.25}
\begin{aligned}
\sum_{n=1}^\iy(\Dl x_n)x_n&=\sum_{n=1}^\iy (x_{n+1}-x_n)x_n=\sum_{n=1}^\iy(x_nx_{n+1}-x_n^2)\\
&=\sum_{n=1}^\iy\left[\frac{x_{n+1}^2-x_n^2}{2}-\frac{1}{2}(x_{n+1}-x_n)^2\right]=-\frac{x_1^2}{2}-\frac{1}{2}\sum_{n=1}
^\iy(\Dl x_n)^2.
\end{aligned}
\end{equation}

Thus, should the integral \eqref{1.22} and the series \eqref{1.23} (with the above-mentioned requirements to $x(t)$)
converge (or diverge) together, ``Theorem~1.5" would be true in the part $\alpha$)--$\theta$).
But this asmumption is wrong: for $x(t)=\frac{\cos\pi t}{\sqrt t},$ the integral \eqref{1.22} converges and the
series \eqref{1.24}, as one can easily see, diverges.
Thus the analogy between conditions for solvability of the Hartman-Wintner problems for differential and difference
equations is not complete although the statements of these problems are completely analogous.

\section{Preliminaries}\label{Preliminaries}

To prove the result from \secref{introduction}, we only need some generalities from the theory of
non-oscillating difference equations of order 2.
See \cite{1} for  a comprehensive exposition.
For the reader's convenience,  below we present a standard summary of all the needed facts.
Recall that equation \eqref{1.2} is said to be non-oscillating at infinity if all its solutions do not oscillate,
i.e., do not change sign beginning from a certain number.
Therefore, for any fixed FSS $\{u_n,v_n\}_{n=0}^\iy$ of equation \eqref{1.2} one can  assume that the
inequalities $u_n>0,\ v_n>0$ hold for all $n\ge0.$
Indeed, otherwise we would just have to change numeration in \eqref{1.2} and for each solution find the corresponding
constant factor $\tau=\pm1.$
The last line of relations in \eqref{1.3} completely characterizes the principal  $\{u_n\}_{n=0}^\iy$ and non-principal
  $\{v_n\}_{n=0}^\iy$ solutions of \eqref{1.2}.
Here the principal solution $\{u_n\}_{n=0}^\iy$ is determined uniquely
 up to a costant factor. Let $\{v_n\}_{n=0}^\iy$ be any non-principal solution.
Then the equality
\begin{equation}\label{2.1}
u_n=v_n\sum_{k=n}^\iy\frac{1}{r_kv_kv_{k+1}},\quad n=0,1,2,\dots
\end{equation}
determines the principal solution.
The following inequalities are immediate consequences of \eqref{1.2}, and we will use them repeatedly:
\begin{equation}\label{2.2}
\frac{u_n}{v_n}>\frac{u_{n+1}}{v_{n+1}},\quad 1>\frac{u_{n+1}}{u_n}\ \frac{v_n}{v_{n+1}},\quad
r_nu_nv_{n+1}>1,\quad n=0,1,2,\dots
\end{equation}

\begin{lem}\label{lem2.1}
The following equalities hold:
\begin{equation}\label{2.3}
\sum_{n=0}^\iy\frac{1}{r_nu_nv_{n+1}}=\iy,\qquad\qquad \sum_{n=0}^\iy\frac{1}{r_nv_nu_{n+1}}=\iy.
\end{equation}
Moreover, if $v_{n+1}\ge v_n$ for all $n\ge 0,$ then
\begin{equation}\label{2.4}
\sum_{n=0}^\iy\frac{1}{r_nu_nv_n}=\iy.
\end{equation}
\end{lem}

\begin{remark}\label{rem2.1} We do not use equalities \eqref{2.3}--\eqref{2.4} in the proofs.
We present them here because they enable one to evaluate the requirements from the statements in
\secref{introduction} more precisely.
\end{remark}

 \begin{proof}[Proof of \lemref{lem2.1}]
 For $n\ge0$ from \eqref{2.1} it follows that (see also \eqref{1.3})
\begin{align*}R_n&\doe \sum_{s=n}^\iy \frac{1}{r_su_sv_{s+1}}=\sum_{s=n}^\iy
\frac{1}{r_sv_sv_{s+1}}\left(\sum_{k=s}^\iy\frac{1}{r_kv_kv_{k+1}}\right)^{-1}\\
&\ge \sum_{s=n}^\iy
\frac{1}{r_sv_sv_{s+1}}\left(\sum_{k=n}^\iy\frac{1}{r_kv_kv_{k+1}}\right)^{-1}=\left(\sum_{k=n}^\iy\frac{1}{r_kv_kv_{k+1}}
\right)^{-1}\cdot\left(\sum_{s=n}^\iy\frac{1}{r_sv_sv_{s+1}}\right)=1.
\end{align*}
Since the remainder $R_n$ of the first series in \eqref{2.3} does not converge to zero as $n\to\iy,$ the series
diverges.
Together with inequalities \eqref{2.2} and the comparison theorem for series, this implies that the second series in
\eqref{2.3} and the series \eqref{2.4} diverge.
\end{proof}

\begin{remark}\label{rem2.2}
Throughout Sections \ref{ProblemEquiv}--\ref{examples} below, in the whole proof the letter $\tau$ stands for absolute
positive constants which are not essential for exposition and may differ even within a single chain of computations.

\end{remark}

\newpage

\section{A Problem Equivalent to the Hartman-Wintner Problem}\label{ProblemEquiv}

Our  study of problem \eqref{1.4}--\eqref{1.6} relies upon the following assertion.

\begin{lem}\label{lem3.1}
Problem  \eqref{1.4}--\eqref{1.6} is solvable if and only if the equation
\begin{equation}\label{3.1}
\Delta(r_{n-1}u_{n-1}u_n\Delta\beta_{n-1})=\sigma_nu_n^2\beta_n,\quad n= 1,2,\dots
\end{equation}
has a solution $\{\beta_n\}_{n=0}^\iy$ such that
\begin{equation}\label{3.2}
\lim_{n\to\iy}\beta_n=1,\quad \lim_{n\to\iy}r_nv_nu_{n+1}\Delta\beta_n=0.
\end{equation}
\end{lem}

\renewcommand{\qedsymbol}{}
\begin{proof}[Proof of \lemref{lem3.1}] \ {\it Necessity}.
Suppose that problem \eqref{1.4}--\eqref{1.6} is solvable.
Set
$$\beta_n=\frac{\tilde u_n}{u_n},\quad n=0,1,2,\dots$$
and show that $\{\beta_n\}_{n=0}^\iy$ is a solution of the problem \eqref{3.1}--\eqref{3.2}.
From \eqref{1.4} it follows that $\beta_n\to1$ as $n\to\iy.$
{}From \eqref{1.5} it follows that
\begin{align*}
\frac{\tilde u_{n+1}}{u_{n+1}}&=\frac{\tilde u_n}{u_n}+\frac{\ve_n}{r_nv_nu_{n+1}}\ \frac{\tilde
u_n}{u_n},\quad\lim_{n\to\iy}\ve_n=0\quad \Rightarrow\quad
\Delta\beta_n=\beta_{n+1}-\beta_n=\frac{\ve_n\beta_n}{r_nv_nu_{n+1}}\\ &\quad \Rightarrow
\lim_{n\to\iy}r_nv_nu_{n+1}\Delta\beta_n=\lim_{n\to\iy}\ve_n\beta_n=\lim_{n\to\iy}\ve_n\cdot\lim_{n\to\iy}\beta_n=0.
\end{align*}

Thus, relations \eqref{3.2} are proved.
Furthermore, the equalities
$$\Delta(r_{n-1}\Delta\tilde u_{n-1})=(q_n+\sigma_n)\tilde u_n,\quad \Delta(r_{n-1}\Delta u_{n-1})=q_nu_n,\quad
n=0,1,2,\dots$$
imply equality \eqref{3.1}:
\begin{align*}
  \sigma_nu_n^2\beta_n& =\sigma_nu_n\tilde u_n=u_n\Delta(r_{n-1}\Delta\tilde u_{n-1})-\tilde
 u_n\Delta(r_{n-1}\Delta u_{n-1})\\
 &=r_n(u_n\tilde u_{n+1}-\tilde u_n u_{n+1})+r_{n-1}(u_n\tilde u_{n-1}-\tilde u_nu_{n-1})
  \\
 &=r_nu_nu_{n+1}(\beta_{n+1}-\beta_n)-r_{n-1}u_{n-1}u_n(\beta_n-\beta_{n-1})=\Delta(r_{n-1}u_{n-1}u_n\Delta
 \beta_{n-1}).
\end{align*}
\end{proof}

\renewcommand{\qedsymbol}{\openbox}
\begin{proof}[Proof of \lemref{lem3.1}] \ {\it Sufficiency}.
Suppose that problem \eqref{3.1}--\eqref{3.2} is solvable, and let $\{\beta_n\}_{n=0}^\iy$ be its solution.
Set
\begin{equation}\label{3.3}
\tilde u_n=\beta_nu_n,\quad n=0,1,2,\dots
\end{equation}

Then
\begin{align*}
\Delta(r_{n-1}&u_{n-1}u_n\Delta\beta_{n-1}) =\Delta\left[r_{n-1}u_{n-1}u_n\left(\frac{\tilde
u_n}{u_n}-\frac{\tilde u_{n-1}}{u_{n-1}}\right)\right]=\Delta\left[r_{n-1}(\tilde u_nu_{n-1}-\tilde
u_{n-1}u_n)\right]\\
&=r_n(\tilde u_{n+1}u_n-\tilde u_nu_{n+1})-r_{n-1}(\tilde u_nu_{n-1}-\tilde u_{n-1}u_n)\\
&=u_n[r_n(\tilde u_{n+1}-\tilde u_n)-r_{n-1}(\tilde u_n-\tilde u_{n-1})]-\tilde
u_n[r_n(u_{n+1}-u_n)-r_{n-1}(u_n-u_{n-1})]\\
&=u_n\Delta(r_{n-1}\Delta\tilde u_{n-1})-\tilde u_n\Delta(r_{n-1}\Delta u_{n-1})=u_n\Delta(r_{n-1}\Delta\tilde
u_{n-1})-q_nu_n\tilde u_n\\
&=\sigma_nu_n^2\frac{\tilde u_n}{u_n}=\sigma_nu_n\tilde u_n.
\end{align*}

Since $u_n\ne0$ for $n\ge0,$ the last equality implies
$$\Delta(r_{n-1}\Delta\tilde u_{n-1})=(q_n+\sigma_n)\tilde u_n,\quad n=0,1,2,\dots$$
Here
\begin{equation}\label{3.4}
\lim_{n\to\iy}\frac{\tilde u_n}{u_n}=\lim_{n\to\iy}\beta_n=1.
\end{equation}
Let $|\beta_n|\ge\frac{1}{2}$  for $n\ge n_0\gg 1.$
Then (see \eqref{3.3}) we have $\tilde u_n\ne0$ for $n\ge n_0.$
Set (see \eqref{2.1})
\begin{equation}\label{3.5}
\tilde v_{n+1}=\tilde u_{n+1}\sum_{k=n_0}^n\frac{1}{r_k\tilde u_k\tilde u_{k+1}}\quad\text{for}\quad n\ge n_0.
\end{equation}
By \eqref{3.5} we get
\begin{equation}\label{3.6}
\frac{\tilde v_{n+1}}{\tilde u_{n+1}}-\frac{\tilde v_n}{\tilde u_n}=\frac{1}{r_n\tilde u_n\tilde
u_{n+1}}\quad \Rightarrow\quad r_n(\tilde v_{n+1}\tilde u_n-\tilde u_{n+1}\tilde v_n)=1,\quad\text{for}\quad n\ge n_0.
\end{equation}

{}From \eqref{3.6} we obtain the following chain of equalities:
\begin{align*}
 &\tilde u_{n+1}(r_{n+1}\Delta\tilde v_{n+1})-\tilde u_n(r_n\Delta\tilde v_n)=\tilde v_{n+1}(r_{n+1}\Delta\tilde
 u_{n+1})-\tilde v_n(r_n\Delta \tilde u_n)\Rightarrow\\
 &\tilde u_{n+1}\Delta(r_n\Delta\tilde v_n)+r_n\Delta\tilde u_n\Delta\tilde v_n=\tilde v_{n+1}\Delta(r_n\Delta\tilde
 u_n)+r_n\Delta\tilde u_n\Delta\tilde v_n\Rightarrow\\
&\tilde u_{n+1}\Delta(r_n\Delta\tilde v_n)=\tilde v_{n+1}\Delta(r_n\Delta\tilde u_n)=\tilde
v_{n+1}(q_{n+1}+\sigma_{n+1})\tilde u_{n+1}\Rightarrow\\
&\Delta(r_n\Delta\tilde v_n)=(q_{n+1}+\sigma_{n+1})\tilde v_{n+1},\quad n\ge n_0.
\end{align*}
Thus $\{\tilde u_n,\tilde v_n\}_{n=n_0}^\iy$ is a FSS of equation \eqref{1.1}.
The following obvious equality is a consequence of \eqref{1.3} (see \cite{1}):
\begin{equation}\label{3.7}
v_{n+1}=\tau u_{n+1}+u_{n+1}\sum_{k=n_0}^n\frac{1}{r_ku_ku_{k+1}},\quad n\ge n_0.
\end{equation}
Indeed, from \eqref{1.3}  we get
\begin{equation}
\left.\begin{array}{ll}
&\frac{v_{n+1}}{u_{n+1}}=\frac{v_n}{u_n}+\frac{1}{r_nu_nu_{n+1}}\nonumber\\
&\frac{v_n}{u_n}=\frac{v_{n-1}}{u_{n-1}}+\frac{1}{r_{n-1}u_{n-1}u_n}\nonumber\\
&\hdots\hdots\hdots\hdots\hdots\hdots\hdots\hdots\nonumber\\
&\frac{v_{n_0+1}}{u_{n_0+1}}=\frac{v_{n_0}}{u_{n_0}}+\frac{1}{r_{n_0}u_{n_0}u_{n_0+1}}
\end{array}\right\}
\Rightarrow
\frac{v_{n+1}}{u_{n+1}}=\frac{v_{n_0}}{u_{n_0}}+\sum_{k=n_0}^n\frac{1}{r_ku_ku_{k+1}}\Rightarrow
\eqref{3.7}
\end{equation}

In the following relations, we use \eqref{3.5}, \eqref{3.7}, \eqref{3.4} and Stolz's theorem (see \cite[vol. I]{4}):
\begin{equation}
\begin{aligned}\label{3.8}
\lim_{n\to\iy}\frac{\tilde v_{n+1}}{v_{n+1}}&=\lim_{n\to\iy}\frac{\tilde
u_{n+1}\sum\limits_{k=n_0}^n\frac{1}{r_k\tilde u_k\tilde u_{k+1}}}{\tau
u_{n+1}+u_{n+1}\sum\limits_{k=n_0}^n\frac{1}{r_ku_ku_{k+1}}}\\
&=\lim_{n\to\iy}\frac{\tilde u_{n+1}}{u_{n+1}}\cdot\lim_{n\to\iy}\frac{
 \sum\limits_{k=n_0}^n\frac{1}{r_k\tilde u_k\tilde u_{k+1}}}{\tau
+\sum\limits_{k=n_0}^n\frac{1}{r_ku_ku_{k+1}}} =\lim_{n\to\iy}\frac{u_nu_{n+1}}{\tilde u_n\tilde u_{n+1}}=1.
\end{aligned}
\end{equation}

Let us now check equality \eqref{1.5}:
\begin{align*}
r_nu_nv_n\left[\frac{\tilde u_{n+1}}{\tilde u_n}-\frac{u_{n+1}}{u_n}\right]&=r_nu_nv_n\frac{u_{n+1}}{\tilde
u_n}\left[\frac{\tilde u_{n+1}}{u_{n+1}}-\frac{\tilde u_n}{u_n}\right]\\
&=(r_nv_nu_{n+1}\Delta\beta_n)\frac{1}{\beta_n}\to 0\quad\text{as}\quad n\to\iy.
\end{align*}

It remains to prove equality \eqref{1.6}.
Below we use \eqref{3.6}, \eqref{1.3}, \eqref{1.5}, \eqref{3.4}, and \eqref{3.8}:
\begin{align*}
\frac{\tilde v_{n+1}}{\tilde v_n}&=\frac{\tilde u_{n+1}}{\tilde u_n}+\frac{1}{r_n\tilde u_n\tilde
v_n}=\frac{\tilde u_{n+1}}{\tilde u_n}+\frac{v_{n+1}u_n-u_{n+1}v_n}{\tilde u_n\tilde v_n}\\
&=\frac{\tilde u_{n+1}}{\tilde u_n}+\frac{u_nv_n}{\tilde u_n\tilde
v_n}\left[\frac{v_{n+1}}{v_n}-\frac{u_{n+1}}{u_n}\right]=\frac{\tilde u_{n+1}}{\tilde
u_n}+\left[\frac{u_nv_n}{\tilde u_n\tilde v_n}-1\right]\left[\frac{v_{n+1}}{v_n}-\frac{u_{n+1}}{u_n}\right]\\
&\ +\frac{v_{n+1}}{v_n}-\frac{u_{n+1}}{u_n}=\frac{v_{n+1}}{v_n}+\left[\frac{u_nv_n}{\tilde u_n\tilde
v_n}-1\right]\frac{1}{r_nu_nv_n}+\left[\frac{\tilde u_{n+1}}{\tilde
u_n}-\frac{u_{n+1}}{u_n}\right]=\frac{v_{n+1}}{v_n}+\frac{\delta_n}{r_nu_nv_n},\\
 \delta_n&=\left[\frac{u_nv_n}{\tilde u_n\tilde v_n}-1\right]+\ve_n,\quad \delta_n\to0\quad\text{as}\quad
n\to\iy.
\end{align*}
Here $\{\ve_n\}_{n=0}^\iy$ is the sequence defined by equality \eqref{1.5}:
$$\frac{\tilde u_{n+1}}{\tilde u_n}=\frac{u_{n+1}}{u_n}+\frac{\ve_n}{r_nu_nv_n},\quad\lim_{n\to\iy}\ve_n=0.$$
\end{proof}

\begin{corl}\label{corl3.1.1}
Problem \eqref{3.1}--\eqref{3.2} is solvable if and only if the equation
\begin{equation}\label{3.9}
r_{n-1}u_{n-1}u_n\Delta\beta_{n-1}=-\sum_{k=n}^\iy\sigma_k u_k^2\beta_k,\quad n=1,2,\dots
\end{equation}
has a solution $\{\beta_n\}_{n=0}^\iy$ such that
\begin{equation}\label{3.10}
\lim_{n\to\iy}\beta_n=1,\quad \lim_{n\to\iy}r_nv_nu_{n+1}\Delta\beta_n=0.
\end{equation}
Moreover, the solutions of the problems \eqref{3.1}--\eqref{3.2} and \eqref{3.9}--\eqref{3.10} coincide  if they
exist.
\end{corl}
\renewcommand{\qedsymbol}{}
\begin{proof}[Proof of \corlref{corl3.1.1}] \ {\it Necessity}.
Suppose that problem \eqref{3.1}--\eqref{3.2} is solvable,  let $\{\beta_n\}_{n=0}^\iy$ be its solution, and
let $m\ge1.$
Then
\begin{align}
\sum_{k=n}^{n+m}\Delta(r_{k-1}u_{k-1}u_k\Delta\beta_{k-1})&=r_{n+m}u_{n+m}u_{n+m+1}\Delta\beta_{n+m}-r_{n-1}u_{n-1}
u_n\Delta\beta_{n-1}\nonumber\\
&=\sum_{k=n}^{n+m}\sigma_ku_k^2\beta_k.\label{3.11}
\end{align}
{}From \eqref{3.2} and \eqref{1.3} it follows that
\begin{equation}\label{3.12}
\lim_{n\to\iy}r_nu_nu_{n+1}\Delta\beta_n=\lim_{n\to\iy}(r_nv_nu_{n+1}\Delta\beta_n)\cdot \frac{u_n}{v_n}=0.
\end{equation}
Therefore, as $m\to\iy,$ \eqref{3.11} implies  \eqref{3.9}, and \eqref{3.2} coincides with \eqref{3.10}.
\end{proof}

\renewcommand{\qedsymbol}{\openbox}
\begin{proof}[Proof of \corlref{corl3.1.1}] \ {\it Sufficiency}.
Suppose that problem \eqref{3.9}--\eqref{3.10} is solvable, and let $\{\beta_n\}_{n=0}^\iy$ be its solution.
Then
\begin{align*}
\Delta(r_{n-1}u_{n-1}u_n\Delta\beta_{n-1})&=r_nu_nu_{n+1}\Delta\beta_n-r_{n-1}u_{n-1}u_n\Delta\beta_{n-1}\\
&=-\sum_{k=n+1}^\iy\sigma_ku_k^2\beta_k+\sum_{k=n}^\iy\sigma_ku_k^2\beta_k=\sigma_nu_n^2\beta_n.
\end{align*}
Moreover, \eqref{3.10} coincides with \eqref{3.2}.
\end{proof}

\section{Proof of   Necessary Conditions for Solvability\\ of the Hartman-Wintner Problem}\label{ProofSuf}

In this section, we prove \thmref{thm1.2}.

\begin{proof}[Proof of \thmref{thm1.2}]
If problem \eqref{1.4}--\eqref{1.6} is solvable, then by \lemref{lem3.1} there exists a solution
$\{\beta_n\}_{n=0}^\iy$ of problem \eqref{3.1}--\eqref{3.2}.
Note that from \eqref{1.3} and \eqref{3.2} it follows that
\begin{equation}\label{4.1}
\lim_{n\to\iy} r_nu_n
v_{n+1}\Delta\beta_n=\lim_{n\to\iy}(1+r_nv_nu_{n+1})\Delta\beta_n=\lim_{n\to\iy}\Delta\beta_n+\lim_{n\to\iy}r_nv_nu_{n+1}
\Delta\beta_n=0.
\end{equation}

Let $|\beta_n|\ge\frac{1}{2}$ for $n\ge n_0$ and $m_2\ge m_1\ge n_0+1.$
Denote
$$\alpha_n=r_nu_nu_{n+1}\Delta\beta_n,\quad \tau=\frac{u_0}{v_0}+1.$$
Then from \eqref{3.1} it follows that
\begin{equation}
\begin{aligned}\label{4.2}
Z(m_1,m_2)&\doe\sum_{n=m_1}^{m_2}\sigma_nu_n^2=\sum_{n=m_1}^{m_2}\frac{\Delta
\alpha_{n-1}}{\beta_n}=\sum_{n=m_1}^{m_2}\left[\frac{\alpha_n}{\beta_n}-\frac{\alpha_{n-1}}{\beta_{n-1}}+\frac
{\alpha_{n-1}\Delta\beta_{n-1}}{\beta_{n-1}\beta_n}\right]\\
&=\left.\frac{\alpha_n}{\beta_n}\right|_{m_1-1}^{m_2}+\sum_{n=m_1}^{m_2}\frac{r_{n-1}u_{n-1}u_n(\Delta\beta_{n-1}
)^2}{\beta_{n-1}\beta_n}.
\end{aligned}
\end{equation}

Let $\ve\in(0,1].$ {}From \eqref{3.2} and \eqref{4.1} we conclude
that there is $n_1>n_0$ such that
\begin{equation}\label{4.3}
\sup_{n\ge  n_1}r_nu_nv_{n+1}|\Delta\beta_n|\le \frac{\ve}{16\tau},\quad \sup_{n\ge
n_1}r_nv_nu_{n+1}|\Delta\beta_n|\le\frac{\ve}{16\tau}.
\end{equation}

Using \eqref{4.2}, \eqref{4.3}, \eqref{2.2}, \eqref{2.1}, we now estimate $Z(m_1,m_2)$ for $m_2\ge m_1\ge
n_1+1:$
\begin{align*}
 |Z(m_1,m_2)|&\le 4\sup_{n\ge n_1}
 r_nu_nu_{n+1}|\Delta\beta_n|
+\sum_{n=m_1}^{m_2}
 \frac{r_{n-1}u_{n-1}u_n|\Delta\beta_{n-1}|^2}{|\beta_{n-1}
 \beta_n|}
  \\
&=4\sup_{n\ge
n_1}\frac{u_n}{v_n}(r_nv_nu_{n+1}|\Delta\beta_n|)+4\sum_{n=m_1}^{m_2}\frac{(r_{n-1}u_{n-1}v_n|\Delta\beta_{n-1}|)(
r_{n-1}v_{n-1}u_n|\Delta\beta_{n-1}|)}{r_{n-1}v_{n-1}v_n}\\
&\le 4\frac{u_0}{v_0}\ \frac{\ve}{16\tau}+4\frac{\ve^2}{256\tau^2}\sum_{n=1}^\iy\frac{1}{r_{n-1}v_{n-1}v_n}\le \frac{\ve}{4}+\frac{\ve^2}{64\tau^2}\
\frac{u_0}{v_0}<\ve.
\end{align*}

This estimate implie that the series $\sigma$ (see \eqref{1.12}) converges in view of Cauchy's criterion.
Hence the sequence $\{C_n\}_{n=0}^\iy$ (see \eqref{1.8}) is well-defined.
It remains to verify \eqref{1.13}.
Denote
\begin{align*}
 &a_n=\sup_{m\ge n-1}\max\{r_{m-1}u_{m-1}v_m|\Delta\beta_{m-1}|,\quad r_{m-1}v_{m-1}u_m|\Delta\beta_{m-1}|\},\quad
 n\ge 1
\\
&b_n=\sup_{m\ge n}\left\{\frac{1}{|\beta_m|},\
  \frac{1}{|\beta_m\beta_{m+1}|}\right\},\quad n\ge 0.
\end{align*}
{}From \eqref{3.2} and \eqref{4.1} we get
\begin{equation}\label{4.4}
\lim_{n\to\iy}a_n=0,\quad \lim_{n\to\iy}b_n=1.
\end{equation}

In the following relations we use the above notation and \eqref{1.8}, \eqref{3.12}, and \eqref{2.1}:
\begin{equation}
\begin{aligned}\label{4.5}
|C_n|&=\frac{v_n}{u_n}\left|\sum_{m=n}^\iy\frac{\Delta\alpha_{m-1}}{\beta_m}\right|=\frac{v_n}{u_n}\left|\sum_{
m=n}^\iy\left\{\frac{\alpha_m}{\beta_{m+1}}-\frac{\alpha_{m-1}}{\beta_m}+\frac{\alpha_m\Delta\beta_m}{\beta_m
\beta_{m+1}}\right\}\right|
 \\
&=\frac{v_n}{u_n}\left|\sum_{m=n}^\iy\frac{\alpha_m\Delta\beta_m}{\beta_m\beta_{m+1}}-\frac{\alpha_{n-1}}{\beta_n}
  \right|\le \frac{r_{n-1}u_{n-1}v_n|\Delta\beta_{n-1}|}{\beta_n}\\
&+\frac{v_n}{u_n}\sum_{m=n}^\iy\frac{(r_mu_mv_{m+1}|\Delta\beta_m|)(r_mv_mu_{m+1}|\Delta\beta_m|)}{|\beta_m\beta_
 {m+1}|r_mv_mv_{m+1}}
\\
&\le a_nb_n\left\{1+a_n\frac{v_n}{u_n}\sum_{m=n}^\iy\frac{1}{r_mv_mv_{m+1}}\right\}=a_nb_n(1+a_n).
\end{aligned}
\end{equation}
{}From \eqref{4.5} and \eqref{4.4} we get \eqref{1.13}.
\end{proof}

\section{Auxiliary Assertions}\label{Aux}

In this section we present various technical assertions needed for  the proofs of Theorems \ref{thm1.3} and
\ref{thm1.4} and their corollaries.
Since we use Assertion \ref{assert1}, below we present its proof for the sake of
completeness.

\begin{proof}[Proof of Assertion \ref{assert1}]
Let
\begin{equation}\label{7.1}
R_n=\sum_{k=n}^\iy \sigma_k u_k^2,\quad n=0,1,2,\dots
\end{equation}
(In particular, $R_0=\sigma$ (see \eqref{1.12}).
Below we use notations  \eqref{1.8}--\eqref{1.9}, and relations \eqref{1.3}:
\begin{equation}
\begin{aligned}\label{7.2}
R_n&=\sum_{k=n}^\iy \sigma_ku_kv_k\ \frac{u_k}{v_k}=\sum_{k=n}^\iy(J_k-J_{k+1})\frac{u_k}{v_k}\\
 &=\sum_{k=n}^\iy\left[J_k\
\frac{u_k}{v_k}-J_{k+1}\ \frac{u_{k+1}}{v_{k+1}}+J_{k+1}\left(\frac{u_{k+1}}{v_{k+1}}-\frac{u_k}{v_k}\right)\right]\\
&=\sum_{k=n}^\iy\left\{J_k\ \frac{u_k}{v_k}-J_{k+1}\ \frac{u_{k+1}}{v_{k+1}}-\frac{J_{k+1}}{r_kv_kv_{k+1}}\right\}=J_n\
\frac{u_n}{v_n}-\sum_{k=n}^\iy\frac{J_{k+1}}{r_kv_kv_{k+1}}.
\end{aligned}
\end{equation}

The series in the right-hand side of \eqref{7.2} absolutely converges since by \eqref{2.1} we have
\begin{equation}\label{7.3}
\sum_{k=n}^\iy \frac{|J_{k+1}|}{r_kv_kv_{k+1}}\le A_n\sum_{k=n}^\iy \frac{1}{r_kv_kv_{k+1}}=A_n\  \frac{u_n}{v_n},\quad n\ge
0.
\end{equation}
Hence the series $R_n,$\ $n\ge0$ converges (at least, conditionally), and therefore by \eqref{1.8} the sequence
$\{C_n\}_{n=0}^\iy$ is well-defined.
{}From \eqref{7.2}--\eqref{7.3} we get
$$|C_n|=|\frac{v_n}{u_n}\ R_n|=\left|J_n-\frac{v_n}{u_n}\sum_{k=n}^\iy\frac{J_{k+1}}{r_kv_kv_{k+1}}\right|\le 2A_n,\quad
n\ge0.$$
\end{proof}

\begin{lem}\label{lem7.1a}
Suppose that the series $J$ (see \eqref{1.7}) converges (at least, conditionally).
Then the series
\begin{equation}\label{7.4a}
H_n\doe\sum_{k=n}^\iy\frac{C_{k+1}}{r_ku_kv_{k+1}},\quad n=0,1,2,\dots
\end{equation}
also converges (at least, conditionally), and the following equality holds:
\begin{equation}\label{7.5a}
H_n=\frac{v_n}{u_n}\sum_{k=n}^\iy\frac{J_{k+1}}{r_kv_kv_{k+1}},\quad n=0,1,2,\dots
\end{equation}
\end{lem}

\begin{proof}
According to \eqref{1.8}, we get (see \eqref{7.1} and \eqref{1.3}):
\begin{align*}
\Delta C_n&=C_{n+1}-C_n=\frac{v_{n+1}}{u_{n+1}}\ R_{n+1}-\frac{v_n}{u_n}\
R_n=\left(\frac{v_{n+1}}{u_{n+1}}-\frac{v_n}{u_n}\right)R_{n+1}+\frac{v_n}{u_n}(R_{n+1} -R_n)\\
&=\frac{1}{r_nu_nu_{n+1}}\ R_{n+1}-\sigma_nu_nv_n=\frac{C_{n+1}}{r_nu_nv_{n+1}}-\sigma_nu_nv_n.
\end{align*}
Since the series $J$ converges,   \eqref{1.9} implies
\begin{equation}\label{7.7a}
-C_n=\sum_{k=n}^\iy\Delta C_k=\sum_{k=n}^\iy \frac{C_{k+1}}{r_ku_kv_{k+1}}-J_n\quad\Rightarrow \quad J_n-C_n=H_n.
\end{equation}
On the other hand, from  \eqref{7.2} and \eqref{1.8}, we get
\begin{equation}\label{7.8}
J_n-C_n=\frac{v_n}{u_n}\sum_{k=n}^\iy\frac{J_{k+1}}{r_kv_kv_{k+1}}.
\end{equation}
Thus the series $H_n,$\ $n=0,1,2,\dots$ converges because of \eqref{7.7a}; and \eqref{7.5a} follows from \eqref{7.7a} and
\eqref{7.8}.
\end{proof}

\begin{proof}[Proof of \lemref{lem1.1}]
Since the series $J$ (see \eqref{1.7}) converges (at least, conditionally) by Assertion \ref{assert1}, the sequence
$\{C_n\}_{n=0}^\iy$ is well-defined (see \eqref{1.8}). Denote
$$\mu_n=\Re C_n,\quad \eta_n=\Im C_n,\quad n=0,1,2,\dots$$
\begin{equation}
\begin{aligned}\label{8.1}
&L(m_1,m_2)=\sum_{n=m_1}^{m_2}\frac{\Re(J_{n+1}\ov C_{n+1})}{r_nu_nv_{n+1}},\quad
G(m_1,m_2)=\sum_{n=m_1}^{m_2}\frac{|C_{n+1}|^2}{r_nu_nv_{n+1}}\\
&\tilde L(m_1,m_2)=\sum_{n=m_1}^{m_2}\frac{\mu_{n+1}}{r_nu_nv_{n+1}}\sum_{k=n+1}^\iy\frac{\mu_{k+1}}{r_ku_kv_{k+1}},\quad
m_2\ge m_1\ge0\\
&\hat L(m_1,m_2)=\sum_{n=m_1}^{m_2}\frac{\eta_{n+1}}{r_nu_nv_{n+1}}\sum_{k=n+1}^\iy\frac{\eta_{k+1}}{r_ku_kv_{k+1}},\quad
m_2\ge m_1\ge0\\
&\qquad\quad \alpha_n=\sum_{k=n}^\iy\frac{\mu_{k+1}}{r_ku_kv_{k+1}},\quad
\gamma_n=\sum_{k=n}^\iy\frac{\eta_{k+1}}{r_ku_kv_{k+1}},\quad n\ge 0.
\end{aligned}
\end{equation}
By \lemref{lem7.1a}, the sequences $\{\alpha_n\}_{n=0}^\iy,$\ $\{\gamma_n\}_{n=0}^\iy$  and the values $\tilde L(m_1,m_2)$,
$\hat L(m_1,m_2)$ are well-defined.
Moreover, by
\eqref{7.7a} and \eqref{1.13} we have
\begin{equation}\label{8.2}
\lim_{n\to\iy}\alpha_n=\lim_{n\to\iy}\gamma_n=0.
\end{equation}

Let us now begin the proof.
From \eqref{7.7a} and \eqref{7.4a}, it follows that
\begin{equation}\label{8.3}
J_{n+1}=C_{n+1}+H_{n+1}=C_{n+1}+\sum_{k=n+1}^\iy\frac{C_{k+1}}{r_ku_kv_{k+1}},\quad n\ge0.
\end{equation}
In turn, from \eqref{8.3} we get for $n\ge 0:$
\begin{equation}\label{8.4}
\frac{\Re(J_{n+1}\ov
C_{n+1})}{r_nu_nv_{n+1}}=\frac{|C_{n+1}|^2}{r_nu_nv_{n+1}}+\frac{1}{r_nu_nv_{n+1}}\sum_{k=n+1}^\iy\frac{\Re(\ov
C_{n+1}C_{k+1})}{r_ku_kv_{k+1}}.
\end{equation}
The last equality implies for $m_2\ge m_1\ge0:$
\begin{equation}
\begin{aligned}\label{8.5}
L(m_1,m_2)&=G(m_1,m_2)+\sum_{n=m_1}^{m_2}\frac{1}{r_nu_nv_{n+1}}\sum_{k=n+1}^\iy\frac{\Re(\ov
C_{n+1}C_{k+1})}{r_ku_kv_{k+1}}\\
&=G(m_1,m_2)+\tilde L(m_1,m_2)+\hat L(m_1,m_2).
\end{aligned}
\end{equation}

The expression  $\tilde L(m_1,m_2)$ can be  written in another way (see \eqref{1.25}):
\begin{equation}
\begin{aligned}\label{8.6}
\tilde
L(m_1,m_2)&=\sum_{n=m_1}^{m_2}(\alpha_n-\alpha_{n+1})\alpha_{n+1}=-\frac{1}{2}\sum_{n=m_1}^{m_2}[\alpha_{n+1}^2-\alpha_n^2+(\alpha
_{n+1}-\alpha_n)^2]\\
&=\frac{\alpha_{m_1}^2-\alpha_{m_2+1}^2}{2}-\frac{1}{2}\sum_{n=m_1}^{m_2}\frac{\mu_{n+1}^2}{(r_nu_nv_{n+1})^2}.
\end{aligned}
\end{equation}
Similarly, we get
\begin{equation}\label{8.7}
\hat L(m_1,m_2)=\frac{\gamma_{m_1}^2-\gamma_{m_2+1}^2}{2}-\frac{1}{2}\sum_{n=m_1}^{m_2}\frac{\eta_{n+1}^2}{(r_nu_nv_{n+1})^2}.
\end{equation}

{}From \eqref{8.6}, \eqref{8.7}, \eqref{8.5} and \eqref{7.4a}, we get
\begin{equation}\label{8.8}
L(m_1,m_2)=G(m_1,m_2)+\frac{|H_{m_1}|^2-|H_{m_2+1}|^2}{2}-\frac{1}{2}\sum_{n=m_1}^{m_2}\frac{|C_{n+1}|^2}{(r_nu_nv_{n+1})^2}.
\end{equation}
Denote
$$\varkappa(m_1,m_2)=\frac{|H_{m_1}|^2-|H_{m_2+1}|^2}{2}.$$
{}From \eqref{7.7a} and \eqref{1.13} it follows that
\begin{equation}\label{8.9}
\lim_{m_1\to\iy}\varkappa(m_1,m_2)=0.
\end{equation}
{}From \eqref{8.8} and \eqref{2.2}  it now follows that
$$|L(m_1,m_2))| \ge
\left|G(m_1,m_2)-\frac{1}{2}\sum_{n=m_1}^{m_2}\frac{|C_{n+1}|^2}{(r_nu_nv_{n+1})^2}\right|-|\varkappa(m_1,m_2)|
 \ge \frac{1}{2} G(m_1,m_2)-|\varkappa(m_1,m_2)|.$$
Hence   by \eqref{2.2} and the last inequality, we have
\begin{equation}\label{8.10}
G(m_1,m_2)\le  2|L(m_1,m_2)|+2|\varkappa(m_1,m_2)|
\end{equation}
\begin{equation}\label{8.11}
|L(m_1,m_2)|\le \frac{3}{2} G(m_1,m_2)+|\varkappa(m_1,m_2)|.
\end{equation}
{}From \eqref{8.10}--\eqref{8.11}, we conclude that the series $G$ and $L$  satisfy  (or do not
satisfy) together
 Cauchy's convergence criterion, which  was to be proved.

Let us now prove the second statement of the lemma concerning the series $P$ (see \eqref{1.16}).
Denote
$$\Re\sigma_n=\mu_n,\quad \Im\sigma_n=\eta_n,\quad n\ge 0;$$
$$\alpha_n=\sum_{k=n}^\iy \mu_k u_k^2,\quad \gamma_n=\sum_{k=n}^\iy \eta_k u_k^2,\quad n\ge 0;$$
$$P(m_1,m_2)=\sum_{n=m_1}^{m_2}\Re(\sigma_n\ov C_n) u_nv_n,\quad m_2\ge m_1\ge0;$$
\begin{equation}\label{8.12}
\tilde P(m_1,m_2)=\sum_{n=m_1}^{m_2}\left(\frac{v_n}{u_n}\right)^2\left[\mu_nu_n^2\sum_{k=n}^\iy
\mu_ku_k^2\right],\quad m_2\ge m_1\ge 0;
\end{equation}
$$\hat P(m_1,m_2)=\sum_{n=m_1}^{m_2}\left(\frac{v_n}{u_n}\right)^2\left[\eta_nu_n^2\sum_{k=n}^\iy \eta_ku_k^2\right],\quad
m_2\ge m_1\ge 0;$$
$$\varkappa(m_,m_2)=|C_{m_1}|^2-|C_{m_2}|^2,\quad B(m_1,m_2)=\sum_{n=m_1}^{m_2}(|\sigma_n|u_nv_n)^2,\quad m_2\ge m_1\ge 0.$$

 According to Assertion \ref{assert1}, all values in \eqref{8.12} are well-defined.
{}From  \eqref{1.8} it follows that
\begin{equation}
\begin{aligned}\label{9.9}
\Re(\sigma_nu_nv_n\ov C_n)&=\Re\left(\sigma_nu_nv_n\frac{v_n}{u_n}\sum_{k=n}^\iy
\ov\sigma_ku_k^2\right)=\left(\frac{v_n}{u_n}\right)^2\Re\left(\sigma_nu_n^2\sum_{k=n}^\iy\ov\sigma_ku_k^2\right)\\
&=\left(\frac{v_n}{u_n}\right)^2\left[\mu_nu_n^2\sum_{k=n}^\iy\mu_ku_k^2+\eta_nu_n^2\sum_{k=n}^\iy\eta_ku_k^2\right].
\end{aligned}
\end{equation}
{}From \eqref{9.9} it immediately follows that
\begin{equation}\label{9.10}
P(m_1,m_2)=\tilde P(m_1,m_2)+\hat P(m_1,m_2).
\end{equation}
Note that
\begin{equation}\label{9.11}
\alpha_n=\Re\left(\frac{u_n}{v_n}C_n\right),\quad\gamma_n=\Im\left(\frac{u_n}{v_n}C_n\right),\quad n\ge0,
\end{equation}
and the values $\tilde P(m_1,m_2),$\ $\hat P(m_1,m_2)$ can be writtenin a form different than \eqref{8.12}.
For example (see \eqref{1.25}):
\begin{align*}
2\tilde
P(m_1,m_2)&=2\sum_{n=m_1}^{m_2}\left(\frac{v_n}{u_n}\right)^2[\alpha_n(\alpha_n-\alpha_{n+1})]=\sum_{n=m_1}^{m_2}
\left(\frac{v_n}{u_n}\right)^2[\alpha_n^2-\alpha_{n+1}^2+(\alpha_n-\alpha_{n+1})^2]\\
 &=\sum_{n=m_1}^{m_2}\left[\left(\frac{v_n\alpha_n}{u_n}\right)^2
 -\left(\frac{v_{n+1}\alpha_{n+1}}{u_{n+1}}\right)^2+\alpha_{n+1}
  ^2\left(\frac{v_{n+1}^2}{u_{n+1}^2}-\frac{v_n^2}{u_n^2}\right)+(\mu_n u_nv_n)^2\right]\\
&=(\Re C_{m_1})^2-(\Re
C_{m_2+1})^2+\sum_{n=m_1}^{m_2}\frac{\alpha_{n+1}^2}{r_nu_nu_{n+1}}\left(\frac{v_{n+1}}{u_{n+1}}+\frac{v_n}{u_n}\right)
+\sum_{n=m_1}^{m_2}(\mu_nu_nv_n)^2\\
&=(\Re C_{m_1})^2-(\Re C_{m_2+1})^2+\sum_{n=m_1}^{m_2}\frac{(\Re
C_{n+1})^2}{r_nu_nv_{n+1}}\left(1+\frac{v_nu_{n+1}}{v_{n+1}u_n}\right)+\sum_{n=m_1}^{m_2}(\mu_nu_nv_n)^2.
\end{align*}

Similarly, we get
$$2\hat
P(m_1,m_2)
=(\Im C_{m_1})^2-(\Im C_{m_2+1})^2+\sum_{n=m_1}^{m_2}\frac{(\Im C_{n+1})^2}{r_nu_nv_{n+1}}\left(1+\frac{u_{n+1}}{u_n}\
\frac{v_n}{v_{n+1}}\right)
+\sum_{n=m_1}^{m_2}(\eta_nu_nv_n)^2.$$
Hence
$$2P(m_1,m_2)=\varkappa(m_1,m_2)+\sum_{n=m_1}^{m_2}\frac{|C_{n+1}|^2}{r_nu_nv_{n+1}}\left(1+\frac{v_n}{v_{n+1}}\
\frac{u_{n+1}}{u_n}\right)+B(m_1,m_2).$$
The last equality and \eqref{2.2} imply
\begin{equation}\label{9.12}
|P(m_1,m_2)|\le|\varkappa(m_1,m_2)|+G(m_1,m_2)+B (m_1,m_2);
\end{equation}
\begin{equation}\label{9.13}
2|P(m_1,m_2)|+|\varkappa(m_1,m_2)|\ge G(m_1,m_2)+B(m_1,m_2).
\end{equation}

Since
\begin{equation}\label{9.14}
\lim_{m_1\to\iy}\varkappa(m_1,m_2)=0,
\end{equation}
the assertion of the lemma now follows from \eqref{9.12}, \eqref{9.13}, \eqref{9.14} and Cauchy's convergence criterion.
\end{proof}

 \begin{lem}\label{lem5.2}
Suppose that the series $J$ (see \eqref{1.7}) converges (at least, conditionally).
Then the series $G$ (see \eqref{1.14}) converges if   the series
\begin{equation}\label{5.26}
I=\sum_{n=0}^{\iy}\frac{|J_{n+1}|^2}{r_nu_nv_{n+1}}
\end{equation}
converges.
\end{lem}
\begin{proof}
We use equality \eqref{7.7a} (see also \eqref{7.4a}):
\begin{equation}\label{9.2}
J_{n+1}=C_{n +1}+\sum_{k=n+1}^\iy\frac{C_{k+1}}{r_ku_kv_{k+1}}=C_{n+1}+H_{n+1}.
\end{equation}
{}From \eqref{9.2} we get (using notations \eqref{8.1}):
\begin{equation}\label{9.3}
(\Re J_{n+1})^2=(\Re C_{n+1})^2+2\mu_{n+1}\sum_{k=n+1}^\iy\frac{\mu_{k+1}}{r_ku_kv_{k+1}}+(\Re H_{n+1})^2;
\end{equation}
\begin{equation}\label{9.4}
(\Im J_{n+1})^2=(\Im C_{n+1})^2+2\eta_{n+1}\sum_{k=n+1}^\iy\frac{\eta_{k+1}}{r_ku_kv_{k+1}}+(\Im H_{n+1})^2.
\end{equation}

Let $m_2\ge m_1\ge0.$ Set
$$
I(m_1,m_2)=\sum_{n=m_1}^{m_2}\frac{|J_{n+1}|^2}{r_nu_nv_{n+1}},\quad
H(m_1,m_2)=\sum_{n=m_1}^{m_2}\frac{|H_{n+1}|^2}{r_nu_nv_{n+1}}
$$

{}From \eqref{9.3}--\eqref{9.4} it immediately follows that
$$
I(m_1,m_2)=G(m_1,m_2)+H(m_1,m_2)+\hat L(m_1,m_2)+\tilde L(m_1,m_2)
$$
(see \eqref{8.1}).
{}From the last inequality by \eqref{8.6}, \eqref{8.7}, \eqref{8.8}, and \eqref{8.9}, we get
\begin{equation}\label{9.7}
I(m_1,m_2)=G(m_1,m_2)-\frac{1}{2}\sum_{n=m_1}^{m_2}\frac{|C_{n+1}|^2}{(r_nu_nv_{n+1})^2}+H(m_1,m_2)+\varkappa(m_1,m_2).
\end{equation}
{}From \eqref{9.7}, by \eqref{2.2}, we get
  \begin{align}
I(m_1,m_2)&+|\varkappa(m_1,m_2)|
\ge|G(m_1,m_2)-\frac{1}{2}\sum_{n=m_1}^{m_2}\frac{|C_{n+1}|^2}{(r_nu_nv_{n+1})^2}+H(m_1,m_2)|\nonumber\\
&=G(m_1,m_2)-\frac{1}{2}\sum_{n=m_1}^{m_2}\frac{|C_{n+1}|^2}{(r_nu_nv_{n+1})^2}+H(m_1,m_2)\ge\frac{1}{2}G(m_1,m_2)+H(m_1,m_2)
\nonumber\\
&\ge\frac{1}{2}G(m_1,m_2).\label{9.8}
\end{align}

{}From \eqref{8.9}, \eqref{9.8} and Cauchy's convergence criterion
we now conclude that if the series $I$ converges (see
\eqref{5.26}), then the series $G$ also converges.
\end{proof}

\section{Proof of the Criterion for Equivalence of the Hartman-Wintner Problem to its Restriction}\label{proofcrit}

In this section we prove \thmref{thm1.4}.

\renewcommand{\qedsymbol}{}
\begin{proof}[Proof of \thmref{thm1.4}] \ {\it Necessity}.

We need  the following obvious assertions.
\end{proof}

\begin{lem}\label{lem6.1}
Problem \eqref{1.1}--\eqref{1.$*$} is solvable if and only if problem \eqref{3.9}--\eqref{3.10} has a solution
$\{\beta_n\}_{n=0}^\iy$ such that
\begin{equation}\label{5.1}
\sum_{n=0}^\iy r_nu_nv_{n+1}|\Delta\beta_n|^2<\iy.
\end{equation}
\end{lem}

\renewcommand{\qedsymbol}{}
\begin{proof}[Proof of \lemref{lem6.1}] \ {\it Necessity}.
Suppose that problem \eqref{1.1}--\eqref{1.$*$} is solvable.
Then, clearly,  problem \eqref{1.4}--\eqref{1.6} is also solvable, and the two problems have  a common solution, a FSS of
equation \eqref{1.1}.
Then by \lemref{lem3.1} and its \corref{corl3.1.1}, equation \eqref{3.9} has a solution $\{\beta_n\}_{n=0}^\iy$ with
properties
\eqref{3.10}, and $\tilde u_n=\beta_nu_n,$\ $n=0,1,2,\dots$ (see \eqref{3.3}).
Therefore \eqref{1.$*$} implies
$$\iy>\sum_{n=0}^\iy r_nu_nv_{n+1}\left|\frac{\tilde u_{n+1}}{u_{n+1}}-\frac{\tilde
u_n}{u_n}\right|^2=\sum_{n=0}^\iy r_nu_nv_{n+1}|\beta_{n+1}-\beta_n|^2=\sum_{n=0}^\iy
r_nu_nv_{n+1}|\Delta\beta_n|^2.
$$
\end{proof}

\renewcommand{\qedsymbol}{\openbox}
\begin{proof}[Proof of \lemref{lem6.1}] \ {\it Sufficiency}.
Suppose that problem \eqref{3.9}--\eqref{3.10} has a solution $\{\beta_n\}_{n=0}^\iy$ with property \eqref{5.1}.
Then by \lemref{lem3.1}  and its \corref{corl3.1.1}, problem \eqref{1.4}--\eqref{1.6} is solvable, and $\tilde u_n=\beta_nu_n,$\
$n\ge0.$ Then the series \eqref{5.1} is none other than the series \eqref{1.$*$}.
\end{proof}

\begin{lem}\label{lem6.2}
The Hartman-Wintner problem is solvable if and only if the series $\sigma$ (see \eqref{1.12}) converges (at least,
conditionally), and the equation
 \begin{equation}\label{6.2a}
 r_nu_nu_{n+1}\Delta
 \beta_n=-\beta_{n+1}\sum_{k=n+1}^\iy\sigma_ku_k^2-\sum_{k=n+1}^\iy\Delta\beta_k\sum_{m=k+1}^\iy\sigma_mu_m^2
 \end{equation}
has a solution $\{\beta_n\}_{n=0}^\iy$ such that
\begin{equation}\label{6.3a}
 \lim_{n\to\iy}\beta_n=1,\qquad \lim_{n\to\iy}r_nv_nu_{n+1}\Delta\beta_n=0.
 \end{equation}
Moreover, the solutions of problem \eqref{3.1}--\eqref{3.2} and \eqref{6.2a}--\eqref{6.3a} coincide if they exist.
In addition, the narrow Hartman-Wintner problem is solvable if and only if the series $\sigma$ converges (at least,
conditionally) and there exists a solution of problem \eqref{6.2a}--\eqref{6.3a} with property \eqref{5.1}
\end{lem}

\renewcommand{\qedsymbol}{}
\begin{proof}[Proof of \lemref{lem6.2}] \ {\it Necessity}.
Suppose that problem \eqref{1.4}--\eqref{1.6} is solvable.
Then by \thmref{thm1.2} the series $\sigma$ converges (at least, conditionally), and by \lemref{lem3.1} and its
\corref{corl3.1.1} there exists a solution $\{\beta_n\}_{n=0}^\iy$ of problem  \eqref{3.9}--\eqref{3.10}. Denote
$$
\alpha_n=\sum_{k=n}^\iy\sigma_ku_k^2,\quad  n=0,1,2,\dots
$$

By \thmref{thm1.2} the sequence $\{\alpha_n\}_{n=0}^\iy$ is well-defined.
Let us transform \eqref{3.9}:
\begin{align}
r_nu_nu_{n+1}\Delta\beta_n&=-\sum_{k=n+1}^\iy\sigma_k
u_k^2\beta_k=-\sum_{k=n+1}^\iy(\alpha_k-\alpha_{k+1})\beta_k\nonumber\\
&=-\sum_{k=n+1}^\iy[\beta_k\alpha_k-\beta_{k+1}\alpha_{k+1}+\alpha_{k+1}\Delta\beta_k]=-\beta_{n+1}\alpha_{n+1}-\sum
_{k=n+1}^\iy(\Delta\beta_k)\alpha_{k+1}\nonumber\\
&=-\beta_{n+1}\sum_{k=n+1}^\iy\sigma_ku_k^2-\sum_{k=n+1}^\iy\Delta\beta_k\left(\sum_{m=k+1}^\iy\sigma_mu_m^2\right).\label{5.3}
\end{align}

The rest of the assertions of the lemma either hold automatically, or follow from \lemref{lem6.1}.
\end{proof}

\renewcommand{\qedsymbol}{\openbox}
\begin{proof}[Proof of \lemref{lem6.2}] \ {\it Sufficiency}.

Since the series $\sigma$ converges and there exists a solutions $\{\beta_n\}_{n=0}^\iy$ of problem
\eqref{6.2a}--\eqref{6.3a}, we have well-defined values
$$r_nu_nu_{n+1}\Delta\beta_n,\qquad \beta_{n+1}\sum_{k=n+1}^\iy\sigma_ku_k^2,\quad n=0,1,2,\dots$$
Then by \eqref{6.2a} we conclude that the series
$$\sum_{k=n+1}^\iy\Delta\beta_k\sum_{m=k+1}^\iy\sigma_mu_m^2$$
converges, and therefore the chain of computations in \eqref{5.3} can be reversed.
Equality \eqref{3.9} is thus proved.
Thus there exists a solution $\{\beta_n\}_{n=0}^\iy$ of problem \eqref{3.9}--\eqref{3.10}, and it coincides with the solution
of problem \eqref{6.2a}--\eqref{6.3a}.
It remains to consider \lemref{lem3.1}, \corref{corl3.1.1} and \lemref{lem6.1}.
\end{proof}

Let us now return to \thmref{thm1.4}.
Suppose that problem \eqref{1.4}--\eqref{1.6} is solvable.
Then by \thmref{thm1.2} the series $\sigma$ converges.
Since the inequality \eqref{1.$*$} also holds, problem \eqref{6.2a}--\eqref{6.3a} is also solvable and its solution
$\{\beta_n\}_{n=0}^\iy$ satisfies inequality \eqref{5.1} (see \lemref{lem6.1} and \lemref{lem6.2}).
Let us rewrite equality \eqref{6.2a} in terms of
$\{C_n\}_{n=0}^\iy$ (see \eqref{1.8}).
We get

\begin{equation}\label{5.4}
r_nu_nv_{n+1}\Delta\beta_n=-C_{n+1}\beta_{n+1}-\frac{v_{n+1}}{u_{n+1}}\sum_{k=n+1}^\iy\Delta\beta_k\frac{u_{k+1}}{v_{
k+1}}C_{k+1},\quad n\ge 0.
\end{equation}
Denote
\begin{equation}
\begin{aligned}\label{5.5}
&S_n=\sum_{k=n}^\iy r_ku_kv_{k+1}|\Delta\beta_k|^2;\quad
F_n=\frac{v_{n+1}}{u_{n+1}}\sum_{k=n+1}^\iy|\Delta\beta_k|\frac{u_{k+1}}{v_{k+1}}|C_{k+1}|,\quad n\ge 0\\
&G_m^{(\theta)}=\sum_{n=m}^{m+\theta}\frac{|C_{n+1}|^2}{r_nu_nv_{n+1}},\ m\ge 0,\ \theta\ge 0;\quad
\xi_n=\sum_{k=n}^\iy \frac{u_{k+1}}{v_{k+1}}\frac{|C_{k+1}|^2}{r_ku_kv_{k+1}},\quad n\ge 0.
\end{aligned}
\end{equation}
Note that the following equality holds:
\begin{equation}\label{5.6}
\lim_{n\to\iy}\frac{v_n}{u_n}\xi_n=0.
\end{equation}
Indeed, from \eqref{2.1} and \eqref{2.2} it follows that
\begin{equation}\label{5.7}
0\le\xi_n=\sum_{k=n}^\iy\frac{u_{k+1}}{v_{k+1}}\ \frac{v_k}{u_k}\ \frac{|C_{k+1}|^2}{r_kv_kv_{k+1}}\le
\frac{u_n}{v_n}\sup_{k\ge n}|C_k|^2.
\end{equation}
{}From \eqref{5.7} and \eqref{1.13} we obtain \eqref{5.6}.

Let us now consider \eqref{5.4}.
Let $|\beta_n|\ge \frac{1}{2}$ for $n\ge n_0.$
For   such $n$ from \eqref{5.4} it follows that
\begin{equation}\label{5.8}
\frac{|C_{n+1}|^2}{8}\le\frac{|C_{n+1}|^2|\beta_{n+1}|^2}{2}\le(r_nu_nv_{n+1}|\Delta\beta_n|)^2+\left[\frac{v_{n+1}}
{u_{n+1}}\sum_{k=n+1}^\iy |\Delta\beta_k|\frac{u_{k+1}}{v_{k+1}}|C_{k+1}|\right]^2.
\end{equation}
Let $m\ge n_0,$\ $\theta\in N.$
{}From \eqref{5.8} we get
\begin{equation}\label{5.9}
\frac{1}{8}G_m^{(\theta)}\le S_m-S_{m+\theta+1}+\sum_{n=m}^{m+\theta}\frac{F_n^2}{r_nu_nv_{n+1}}\le
S_m+\sum_{n=m}^{m+\theta}\frac{F_n^2}{r_nu_nv_{n+1}}.
\end{equation}
Let us estimate the second sum in \eqref{5.9}.
First note that from the definition of $F_n,$ \ $S_n$ and Schwarz's inequality, we obtain
\begin{equation}\label{5.10}
F_n\le
\frac{v_{n+1}}{u_{n+1}}S_{n+1}^{1/2}\left[\sum_{k=n+1}^\iy\left(\frac{u_{k+1}}{v_{k+1}}\right)^2\frac{|C_{k+1}|^2}
{r_ku_kv_{k+1}}\right]^{1/2}.
\end{equation}

In the following relations, we use \eqref{5.10}, \eqref{2.2}, \eqref{1.3} and \eqref{5.7}:
\begin{equation}
\begin{aligned}\label{5.11}
\sum_{n=m}^{m+\theta}&\frac{F_n^2}{r_nu_nv_{n+1}} \le
\sum_{n=m}^{m+\theta}\frac{S_{n+1}}{r_nu_nv_{n+1}}\left(\frac{
v_{n+1}}{u_{n+1}}\right)^2\sum_{k=n+1}^\iy\left(\frac{u_{k+1}}{v_{k+1}}\right)^2\frac{|C_{k+1}|^2}{r_ku_kv_{k+1}}\\
&\le S_m\sum_{n=m}^{m+\theta}\frac{1}{r_nu_nu_{n+1}}\left[\frac{v_{n+1}}{u_{n+1}}\sum_{k=n+1}^\iy\left(\frac{u_{
k+1}}{v_{k+1}}\right)^2\frac{|C_{k+1}|^2}{r_ku_kv_{k+1}}\right]\\
&\le S_m\sum_{n=m}^{m+\theta}\frac{1}{r_nu_nu_{n+1}}\left[\frac{v_{n+1}}{u_{n+1}}\ \frac{u_{n+2}}{v_{n+2}}\sum_
{k=n+1}^\iy \frac{u_{k+1}}{v_{k+1}}\ \frac{|C_{k+1}|^2}{r_ku_kv_{k+1}}\right]\\
&\le S_m\sum_{n=m}^{m+\theta}  \left[\frac{v_{n+1}}{u_{n+1}}-\frac{v_n}{v_n}\right]\xi_{n+1}=S_m\sum_{n=m}^{m+\theta}
\left[\frac{v_{n+1}}{u_{n+1}}\xi_{n+1}-\frac{v_n}{u_n}\xi_n+\frac{v_n}{u_n}(\xi_n-\xi_{n+1}) \right]\\
&=S_m\left[\left.\frac{v_k}{u_k}\xi_k\right|_m^{m+\theta+1}+\sum_{n=m}^{m+\theta}\frac{v_n}{u_n}\
\frac{u_{n+1}}{v_{n+1}}\ \frac{|C_{n+1}|^2}{r_nu_nv_{n+1}}\right]\le S_m\sup_{k\ge
m+\theta}|C_k|^2+S_mG_m^{(\theta)}.
\end{aligned}
\end{equation}

Let us choose a bigger $n_0$ (if necessary) so that $S_m\le 16^{-1}$ for $m\ge n_0.$
Then \eqref{5.11} and \eqref{5.9} imply
\begin{equation}\label{5.12}
G_m^{(\theta)}\le 16 S_m+16 S_m\sup_{k\ge m+\theta}|C_k|^2,\quad m\ge n_0.
\end{equation}
Let us take the limit in \eqref{5.12} (as $\theta\to\iy$).
Then by \eqref{1.13} we have
$$\sum_{n=m}^\iy\frac{|C_{n+1}|^2}{r_nu_nv_{n+1}}=\lim_{\theta\to\iy}G_m^{(\theta)}\le 16S_m<\iy.$$

\renewcommand{\qedsymbol}{\openbox}
\begin{proof}[Proof of \thmref{thm1.2}] \ {\it Sufficiency}.
Suppose that problem \eqref{1.4}--\eqref{1.6} is solvable and
$G<\iy$ (see \eqref{1.14}). Then by \thmref{thm1.2} the series
$\sigma$ (see \eqref{1.2})   converges (at least, conditionally),
the sequence $\{C_n\}_{n=0}^\iy$   is well-defined, and by
\lemref{lem6.2} there exists a solution of $\{\beta_n\}_{n=0}^\iy$
of problem \eqref{6.2a}--\eqref{6.3a} and it coincides with the
solution of problem \eqref{3.1}--\eqref{3.2}. Let us transform
\eqref{5.4}:
\begin{equation}\label{5.13}
\mu_n=r_nv_nu_{n+1}\Delta\beta_n,\quad n=0,1,2,\dots
\end{equation}
We get
\begin{equation}\label{5.14}
\mu_n=-\frac{u_{n+1}}{u_n}\ \frac{v_n}{v_{n+1}}C_{n+1}\beta_{n+1}-\frac{v_n}{u_n}\sum_{k=n+1}^\iy\frac{C_{k+1}
\mu_k}{r_kv_kv_{k+1}},\quad n\ge 0.
\end{equation}
{}From \eqref{5.14}, taking into account \eqref{3.10}, we get for  $n\ge 0$
\begin{equation}\label{5.15}
|\mu_n| =r_nv_ku_{n+1}|\Delta\beta_n|
\le\tau\left[\frac{u_{n+1}}{u_n}\
\frac{v_n}{v_{n+1}}|C_{n+1}|+\frac{v_n}{u_n}\sum_{k=n+1}^\iy\frac{|C_{k+1}|}{r_kv_kv_{k+1}}\right].
\end{equation}
Let us estimate the sum in the right-hand side of \eqref{5.15}.
Below we use Schwarz's  inequality and \eqref{2.1}:
\begin{equation}
\begin{aligned}\label{5.16}
 \frac{v_n}{u_n}\sum_{k=n+1}^\iy\frac{|C_{k+1}|}{r_kv_kv_{k+1}}&\le
 \frac{v_n}{u_n}\left[\sum_{k=n+1}^\iy\frac{|C_{k+1}|^2}{r_kv_kv_{k+1}}\right]^{1/2}\left[\sum_{k=n+1}^\iy\frac{1}
 {r_kv_kv_{k+1}}\right]^{1/2}\\
&=\frac{v_n}{u_n}\sqrt{\frac{u_{n+1}}{v_{n+1}}}\left[\sum_{k=n+1}^\iy\frac{|C_{k+1}|^2}{r_kv_kv_{k+1}}\right]^{1/2}.
\end{aligned}
\end{equation}

{}From \eqref{5.15} and \eqref{5.16} it follows that
$$r_nv_nu_{n+1}|\Delta\beta_n|\le\tau\left\{\frac{u_{n+1}}{u_n}\
\frac{v_n}{v_{n+1}}|C_{n+1}|+\frac{v_n}{u_n}\sqrt{\frac{u_{n+1}}{v_{n+1}}}\left[\sum_{k=n+1}^\iy\frac{|C_{k+1}|^2}{r_kv
_kv_{k+1}}\right]^{1/2}\right\}.$$
Let us multiply the last inequality by $\sqrt{\frac {u_nv_{n+1}}{r_n}}\ \frac{1}{v_nu_{n+1}}.$
We get
$$\sqrt{r_nu_nv_{n+1}}|\Delta\beta_n|\le \tau\left\{\frac{|C_{n+1}|}{\sqrt{r_nu_nv_{n+1}}}+\frac{1}{\sqrt{r_nu_nu_{n+1}}}\left[\sum_{k=n+1}^\iy\frac{|C_{k+1}|^2}{
r_kv_kv_{k+1}}\right]^{1/2}\right\}.$$
Hence
\begin{equation}\label{5.17}
r_nu_nv_{n+1}|\Delta\beta_n|^2\le\tau\left\{\frac{|C_{n+1}|^2}{r_nu_nv_{n+1}}+\frac{1}{r_nu_nu_{n+1}}\sum_{k=n+1}^\iy
\frac{|C_{k+1}|^2}{r_kv_kv_{k+1}}\right\}.
\end{equation}
Let us sum up the inequalities \eqref{5.17} for $n\ge0.$
Then (see \eqref{5.5} and \eqref{1.14})
\begin{equation}\label{5.18}
S\doe
\sum_{n=0}^\iy r_nu_nv_{n+1}|\Delta\beta_n|^2\le\tau\left\{G+\sum_{n=0}^\iy\frac{1}{r_nu_nu_{n+1}}
\sum_{k=n+1}^\iy\frac{|C_{k+1}|^2}{r_kv_kv_{k+1}}\right\}.
\end{equation}
Let us estimate the sum in the right-hand side of \eqref{5.18}.
Denote
\begin{equation}\label{5.19}
W_n=\sum_{k=n}^\iy \frac{|C_{k+1}|^2}{r_kv_kv_{k+1}},\quad n\ge 0.
\end{equation}
The sequence $\{W_n\}_{n=0}^\iy$ is well-defined.
Indeed, as shown above, the sequence $\{C_n\}_{n=0}^\iy$ is defined and it is bounded because of \eqref{1.13}.
Hence, according to \eqref{2.1}, we have
\begin{equation}\label{5.20}
0\le W_n=\sum_{k=n}^\iy \frac{|C_{k+1}|^2}{r_kv_kv_{k+1}}\le\left(\sup_{k\ge
n}|C_{k+1}|^2\right)\sum_{k=n}^\iy\frac{ 1}{r_kv_kv_{k+1}}=\frac{u_n}{v_n}\sup_{k\ge n}|C_{k+1}|^2.
\end{equation}

{}From \eqref{5.20} and \eqref{1.13} we get
\begin{equation}\label{5.21}
\lim_{n\to\iy}\frac{v_n}{u_n}W_n=0.
\end{equation}
In the following relations, we use \eqref{1.3}, \eqref{5.21}:
\begin{equation}
\begin{aligned}\label{5.22}
\sum_{n=0}^\iy\frac{1}{r_nu_nu_{n+1}}&\sum_{k=n+1}^\iy\frac{|C_{k+1}|^2}{r_kv_kv_{k+1}}=\sum_{n=0}^\iy\left(\frac{v_{n+1}
}{u_{n+1}}-\frac{v_n}{u_n}\right)W_{n+1}\\
&=\sum_{n=0}^\iy\left[\frac{v_{n+1}}{u_{n+1}}W_{n+1}-\frac{v_n}{u_n}W_n-\frac{v_n}{u_n}(W_{n+1}-W_n)\right]\\
&=-\frac{v_0}{u_0}W_0+\sum_{n=0}^\iy\frac{v_n}{u_n}\ \frac{|C_{n+1}|^2}{r_nv_nv_{n+1}}=G-\frac{v_0}{u_0}W_0\le  G.
\end{aligned}
\end{equation}
{}From \eqref{5.17} and \eqref{5.22} it follows that
$$S=\sum_{n=0}^\iy r_nu_nv_{n+1}|\Delta\beta_n|^2\le\tau G<\iy.$$
\end{proof}

\section{Double Purpose Lemmas}\label{double}

In this section we present some technical assertions which will be applied in the proof of \thmref{thm1.3}, both in the
``necessary" and the ``sufficient" part.
In the statements below we use assumptions of two types -- condition A) and condition B):
\begin{enumerate}
\item[A)] the narrow Hartman-Wintner problem (problem \eqref{1.4}--\eqref{1.$*$}) is solvable;
\item[B)] the series $G$ and $J$ converge ($J$ at least, conditionally) (see \eqref{1.4}, \eqref{1.7}).
\end{enumerate}

Throughout the sequel by condition A) we mean the collection of all criteria already proved for solvability of problem
\eqref{1.4}--\eqref{1.$*$}.
Thus the expression ``condition A) holds" means that the following requirements are satisfied:
  \begin{enumerate}
\item[A.1)] the series $\sigma$ (see \eqref{1.12} converges (at least, conditionally) (see \thmref{thm1.2});
\item[A.2)] the sequence $\{C_n\}_{n=0}^\iy$ (see \eqref{1.8}) is well-defined, and equality \eqref{1.13} holds (see\linebreak
\thmref{thm1.2});
\item[A.3)] the inequality $G<\iy$ holds (see \eqref{1.14} and \thmref{thm1.4});
\item[A.4)] the series $\sigma$ converges (at least, conditionally), and problem \eqref{6.2a}--\eqref{6.3a} has a solution
$\{\beta_n\}_{n=0}^\iy$ which satisfies inequality \eqref{5.1} and equality \eqref{5.14} (see \lemref{lem6.2}).
\end{enumerate}

Similarly, the expression ``condition B) holds" means that the following requirements are satisfies:
 \begin{enumerate}
\item[B.1)] the series $J$  converges (at least, conditionally);
\item[B.2)] Condition A.1) holds (see Assertion \ref{assert1});
\item[B.3)] Condition A.2) holds and inequalities \eqref{1.9} hold (see Assertion \ref{assert1});
\item[B.4)] Condition A.3) holds (see \thmref{thm1.4}).\end{enumerate}

Let us proceed to the exposition of the results.

Let $n_0$ denote a number and define
\begin{equation}\label{7.1a}
\theta\doe\{\theta_k\}_{k=n_0}^\iy,\qquad \|\theta\|\doe\sup_{n\ge n_0}|\theta_k|.
\end{equation}
In the space of sequences \eqref{7.1a} let us introduce a linear operator $\A$ and a sequence
$g(\theta)=\{g_n(\theta)\}_{n=n_0}^\iy:$
\begin{equation}\label{8.17}
\mathcal A\theta=\{(\mathcal A\theta)_n\}_{n=n_0}^\iy,\quad (\mathcal
A\theta)_n=\frac{v_n}{u_n}\sum_{k=n+1}^\iy\frac{C_{k+1}\theta_k}{r_kv_kv_{k+1}},\quad n\ge n_0;
\end{equation}
\begin{equation}\label{8.18}
g(\theta)=\{g_n(\theta)\}_{n=n_0}^\iy,\quad g_n(\theta)=-\frac{u_{n+1}}{u_n}\ \frac{v_n}{v_{n+1}}C_{n+1}\theta_{n+1},\quad n\ge
n_0.
\end{equation}
Here $\{C_n\}_{n=n_0}^\iy$ is the sequence defined in \eqref{1.8}.

\begin{lem}\label{lem7.1}
Suppose that either one of the conditions \rm{A)} or \rm{B)}  holds. Then:
\begin{equation}\label{7.4}
\|g(\theta)\|\le\left(\sup_{n\ge n_0}|C_n|\right)\|\theta\|;
\end{equation}
\begin{equation}\label{7.5}
|(\A\theta)_n|\le\left(\sup_{k\ge n}|C_k|\right)\left(\sup_{k\ge n+1}|\theta_k|\right)\le \left(\sup_{k\ge
n}|C_k|\right)\|\theta\|;
\end{equation}
\begin{equation}\label{7.6}
\|\A\|\le\sup_{n\ge n_0}|C_n|.
\end{equation}
\end{lem}

\begin{proof}
In the following relations we use the hypotheses of the lemma, \eqref{2.2} and \eqref{2.1}:
\begin{align*}
&\qquad\qquad \|g(\theta)\|=\sup_{n\ge n_0}\frac{u_{n+1}v_n}{u_nv_{n+1}}|C_{n+1}|\ |\theta_{n+1}|\le\left(\sup_{n\ge
n_0}|C_n|\right)\|\theta\|;\\
&|(\A\theta)_n|\le \frac{v_n}{u_n}\sum_{k=n+1}^\iy\frac{|C_{k+1}|\ |\theta_k|}{r_kv_kv_{k+1}}\le \left(\sup_{k\ge
n}|C_k|\right)\left(\sup_{k\ge n+1}|\theta_k|\right)\frac{v_n}{u_n}\sum_{k=n+1}^\iy\frac{1}{r_kv_kv_{k+1}}\\
&=\left(\sup_{k\ge n}|C_k|\right)\left(\sup_{k\ge n+1}|\theta_k|\right)\frac{v_n}{u_n}\ \frac{u_{n+1}}{v_{n+1}}\le
\left(\sup_{k\ge n}|C_k|\right)\left(\sup_{k\ge k+1}|\theta_k|\right)\le \left(\sup_{k\ge n}|C_k|\right)\|\theta\|;\\
&\qquad\qquad \|\A\theta\|=\sup_{k\ge n_0}|(\A\theta)_k|\le\sup_{k\ge n_0}\left(\sup_{m\ge
k}|C_m|\right)\|\theta\|\le\left(\sup_{k\ge k_0}|C_k|\right)\|\theta\|.
\end{align*}
\end{proof}

\begin{lem}\label{lem7.2}
 Suppose that either one of the conditions  \rm{A)} or \rm{B)} holds.
Then there is $n_0\gg 1$ such that for all $n\ge n_0$ and any vector $\theta$ with finite norm \eqref{7.1a}, the following
equality holds:
  \begin{equation}\label{7.7}
    \Bigg(\sum_{k=1}^\iy(-1)^k\A^kg(\theta)\Bigg)_n
=\sum_{k=1}^\iy(-1)^k(\A^kg(\theta))_n.
 \end{equation}
\end{lem}

\begin{proof}
Let $n\ge n_0$, and let $P_n$ be the functional defined by the equality
$$P_n\theta=\theta_n,\qquad \theta=\{\theta_k\}_{k=n_0}^\iy,\qquad \|\theta\|<\iy.$$
Clearly,  $P_n$ is linear.
Obviously, $\|P_n\|=1$ since
$$|P_n\theta|=|\theta_n|\le\|\theta\|$$
and, in addition,
\begin{equation*}
P_n\theta _0=1=\|\theta _0\|,\quad \theta_0 =\{\theta _k^0\}_{k=n_0}^\iy,\quad \theta _k^0=\begin{cases}0,\quad &\text{if}\ k\ne
n\\ 1,\quad &\text{if}\ k=n
\end{cases}
\end{equation*}

Let us now verify that if
$$
\theta ^{(0)}=\sum_{k=1}^\iy\theta^{(k)},\qquad \sum_{k=1}^\iy\|\theta^{(k)}\|<\iy,
$$
then the following equality holds:
\begin{equation}\label{8.28}
P_n\left(\sum_{k=1}^\iy\theta ^{(k)}\right)=\sum_{k=1}^\iy P_n\theta^{(k)}.
\end{equation}
Indeed,  let $m\ge1.$ Then
\begin{align*}
|P_n\theta ^{(0)}-\sum_{k=1}^m P_n\theta^{(k)}|&=\left|P_n\left(\theta^{(0)}-\sum_{k=1}^m\theta^{(k)}\right)\right|\le \|P_n\|\
\|\theta^{(0)}-\sum_{k=1}^m\theta^{(k)}\|\\
&\le\sum_{k=m+1}^\iy\|\theta^{(k)}\|\to0\quad{as}\quad m\to\iy;\ \Rightarrow\eqref{8.28}.
\end{align*}

Furthermore, from the hypotheses of the lemma it follows that there is $n_0\gg 1$ such that $|C_n|\le 2^{-1}$ for all $n\ge
n_0.$
Then $\|\A\|\le 2^{-1}$ (see \eqref{7.6}), and from \eqref{7.7} and \lemref{lem7.1}  it follows that
$$\sum_{k=1}^\iy\|(-1)^k\A^kg(\theta)\|\le\sum_{k=1}^\iy\|\A\|^k\|g(\theta)\|\le \sum_{k=1}^\iy
\frac{1}{2^{k+1}}\|\theta\|=\tau\|\theta\|<\iy.$$
It remains to apply \eqref{8.28}.
\end{proof}

Denote (see \eqref{8.17}--\eqref{8.18}):
\begin{equation}\label{8.41}
f_n(\theta)=\frac{1}{r_nv_nu_{n+1}}\sum_{k=1}^\iy(-1)^k(\mathcal A^kg(\theta))_n,\quad n\ge n_0;
\end{equation}
\begin{equation}\label{8.45}
\Phi_n(\theta)=\sum_{k=n}^\iy |f_n(\theta)|,\qquad\varkappa_n(\theta)=\sum_{s=1}^\iy\sum_{k=n}^\iy\frac{|(\mathcal
A^sg(\theta))_k|}{r_kv_ku_{k+1}}<\iy,\quad n\ge n_0,
\end{equation}
\begin{equation}\label{8.47}
\varkappa_n^{(s)}(\theta)=\sum_{k=n}^\iy\frac{|(\mathcal A^sg(\theta))_k|}{r_kv_ku_{k+1}},\quad s\ge 1,\ n\ge
n_0.
\end{equation}

\begin{lem}\label{lem7.3}
Suppose that one of the conditions \rm{A)}  or \rm{B)} holds.
Then there is $n_0\gg 1$ such that for all $n\ge n_0$ and any sequence $\theta$ with finite norm \eqref{7.1a}, the following
inequalities hold:
\begin{equation}\label{7.12}
\varkappa_n^{(s)}(\theta)\le G_n\left(\sup_{k\ge n}|C_k|\right)^{s-1}\|\theta\|,\quad s\ge 1,\ n\ge n_0,
\end{equation}
\begin{equation}\label{7.13}
\varkappa_n(\theta)\le 2G_n \|\theta\|,\quad n\ge n_0,
\end{equation}
\begin{equation}\label{7.14}
\Phi_n(\theta)\le 2G_n \|\theta\|,\quad n\ge n_0.
\end{equation}
Here
\begin{equation}\label{8.50}
G_n \doe \sum_{k=n}^\iy\frac{|C_{k+1}|^2}{r_ku_kv_{k+1}},\quad n\ge 0.
\end{equation}
\end{lem}

\begin{proof}
We first prove inequality \eqref {7.12} for $s=1.$ {}From the hypothesis of the
lemma, we conclude that the sequences $\{G_n\}_{n=0}^\iy,$\ $\{C_n\}_{n=0}^\iy$
are well-defined. In the following relations, we consecutively use
\eqref{8.47}, \eqref{8.17}, \eqref{8.18}, \eqref{5.5}, \eqref{1.3}, \eqref{5.6}
and \eqref{2.2}:
\begin{align*}
\varkappa_n^{(1)}(\theta)&=\sum_{k=n}^\iy\frac{|(\mathcal A
g(\theta))_k|}{r_kv_ku_{k+1}}=\sum_{k=n}^\iy\frac{1}{r_kv_ku_{k+1}}\left|\frac{v_k}{u_k}\sum_{m=k+1}^\iy\frac{C_{m+1}g_m(\theta)
}{r_mv_mv_{m+1}}
\right|\\
&\le\sum_{k=n}^\iy\frac{1}{r_ku_ku_{k+1}}\sum_{m=k+1}^\iy\frac{|C_{m+1}|^2}{r_mu_mv_{m+1}}\
\frac{u_{m+1}}{v_{m+1}}\|\theta\|=\sum_{k=n}^\iy\frac{\xi_{k+1}}{r_ku_ku_{k+1}}\|\theta\|
\\
&=\sum_{k=n}^\iy\left(\frac{v_{k+1}}{u_{k+1}}-\frac{v_k}{u_k}\right)\xi_{k+1}\|\theta\|=\sum_{k=n}^\iy\left[\frac{v_{k+1}
 \xi_{k+1}}{u_{k+1}}-\frac{v_k\xi_k}{u_k}+\frac{v_k}{u_k}(\xi_k-\xi_{k+1})\right]\|\theta\|\\
&=\left[-\frac{v_n}{u_n}\xi_n+\sum_{k=n}^\iy \frac{v_k}{u_k}\ \frac{u_{k+1}}{v_{k+1}}\
\frac{|C_{k+1}|^2}{r_ku_kv_{k+1}}\right]\|\theta\|\le G_n\|\theta\|.
\end{align*}

Let us consider the case $s\ge 2.$
Denote (see \eqref{8.17}):
\begin{equation}
\begin{aligned}\label{8.51}
&\mathcal A^{s-1}g(\theta) =h^{(s)}(\theta)=\{h_n^{(s)}(\theta)\}_{n=n_0}^\iy\quad\Rightarrow\quad h_n^{(s)}=(\mathcal
A^{s-1}g(\theta))_n,\quad n\ge n_0,\quad s\ge 2,\ \\
&\qquad\qquad\qquad\qquad\gamma_k^{(s)}(\theta)=\sum_{m=k}^\iy \frac{|C_{m+1}|\ |h_m^{(s)}|(\theta)}{r_mv_mv_{m+1}},\quad s\ge
2,\ k\ge n_0.
\end{aligned}
\end{equation}

Let us verify that for $\|\theta\|<\iy$ we have the equality
\begin{equation}\label{8.52}
\lim_{k\to\iy}\frac{v_k}{u_k}\gamma_{k+1}^{(s)}(\theta)=0.
\end{equation}

Indeed, by \eqref{1.13}, \eqref{2.1} and \eqref{2.2}, we get
\begin{equation}
\begin{aligned}\label{7.18}
\frac{v_k}{u_k}&\gamma_{k+1}^{(s)}(\theta) \le \frac{v_k}{u_k}\Big(\sup_{m\ge k+1}|C_m|\Big)\Big(\sup_{m\ge
k+1}(|h_m^{(s)})(\theta)|\Big)\sum_{m=k+1}^\iy\frac{1}{r_mv_m v_{m+1}}\\
&=\frac{u_{k+1}}{v_{k+1}}\ \frac{v_k}{u_k}\Big(\sup_{m\ge k+1}|C_m|\Big)\sup_{m\ge k+1}(|h_m^{(s)}(\theta)|)\le
\Big(\sup_{m\ge k+1}|C_m|\Big)\Big(\sup_{m\ge k+1}|h_m^{(s)}(\theta)|\Big).
\end{aligned}
\end{equation}
Furthermore, from \eqref{8.51},   \lemref{lem7.1}, \eqref{8.18}, \eqref{1.13} and \eqref{2.2} it follows that
$$
\sup_{m\ge k+1}(|h_m^{(s)}(\theta)|)\le \|h^{(s)}(\theta)\|=\|\mathcal A^{s-1}g(\theta)\|\le \|\mathcal A\|^{s-1}\|g(\theta)\|
\le \tau\|\theta\|.$$

The last estimate and \eqref{7.18} imply \eqref{8.52}.
In the following relations,  we consecutively use \eqref{8.47}, \eqref{8.51}, \eqref{8.17},
\eqref{1.3},
\eqref{8.51},
\eqref{8.52},
  and \eqref{2.2}:
\begin{align*}
\varkappa_n^{(s)}(\theta)&=\sum_{k=n}^\iy\frac{|(\mathcal A^sg(\theta))_k|}{r_kv_ku_{k+1}}=\sum_{k=n}^\iy \frac{|(\mathcal
Ah^{(s)}(\theta))_k|}{r_kv_ku_{k+1}}=\sum_{k=n}^\iy\frac{1}{r_kv_ku_{k+1}}\left|\frac{v_k}{u_k}\sum_{m=k+1}^\iy\frac{C_{m+1}h_m^{
(s)}(\theta)}{r_mv_mv_{m+1}}\right|\\
& \le\sum_{k=n}^\iy\frac{1}{r_ku_ku_{k+1}}\sum_{m=k+1}^\iy\frac{|C_{m+1}|\
|h_{m}^{(s)}(\theta)|}{r_mv_mv_{m+1}}=\sum_{k=n}^\iy\left(\frac{v_{k+1}}{u_{k+1}}-\frac{v_k}{u_k}\right)\gamma_{k+1}^{(s)}(\theta)\\
&=\sum_{k=n}^\iy\left[\frac{v_{k+1}}{u_{k+1}}\gamma_{k+2}^{(s)}(\theta)-\frac{v_k}{u_k}\gamma_{k+1}^{(s)}(\theta)+\frac{v_{k+1}}{u_{k+1}}(\gamma
_{k+1}^{(s)}(\theta)-\gamma_{k+2}^{(s)}(\theta))\right]\\
&=-\frac{v_n}{u_n}\gamma_{n+1}^{(s)}(\theta)+\sum_{k=n}^\iy\frac{v_{k+1}}{u_{k+1}}\ \frac{|C_{k+2}|\
|h_{k+1}^{(s)}(\theta)|}{r_{k+1}v_{k+1}v_{k+2}}\le\sum_{ k=n+1}^\iy\frac{|C_{k+1}|\ |h_k^{(s)}(\theta)|}{r_ku_kv_{k+1}}\\
&\le \Big(\sup_{k\ge n}|C_k|\Big)\sum_{k=n+1}^\iy \frac{|(\mathcal
A^{s-1}g(\theta))_k|}{r_kv_ku_{k+1}}=\Big(\sup_{k\ge n}|C_k|\Big)\varkappa_{n+1}^{(s-1)}(\theta).
\end{align*}

We thus get the estimate
\begin{equation}\label{7.19}
\varkappa_n^{(s)}(\theta)\le\left(\sup_{k\ge n}|C_k|\right)\varkappa_{n+1}^{(s-1)}(\theta),\quad s\ge 2,\ n\ge n_0.
\end{equation}

To finish the proof of \eqref{7.12}, we consecutively use inequality \eqref{7.19} and the obtained estimate for
$\varkappa_n^{(1)}(\theta):$
  \begin{align*}
  \varkappa_n^{(s)}(\theta)&\le \left(\sup_{k\ge
 n}|C_k|\right)\varkappa_{n+1}^{(s-1)}(\theta)\le\dots\le\left(\sup_{k\ge n}|C_k|\right)^{s-1}\varkappa_{n+s}^{(1)}(\theta)
  \\
 &\le \left(\sup_{k\ge n}|C_k|\right)^{s-1}G_n\|\theta\|.
\end{align*}

{}From the hypothesis of the lemma, it follows that there is $n_0\gg1$ such that $\sup\limits_{k\ge n_0}|C_k|\le 2^{-1}.$
Then \eqref{8.45}, \eqref{8.47} and \eqref{7.12} imply \eqref{7.13}:
  \begin{align*}
  \varkappa_n(\theta)&=\sum_{s=1}^\iy \varkappa_n^{(s)}(\theta)\le G_n\cdot\|\theta\|\sum_{s=1}^\iy\left(\sup_{k\ge
  n_0}|C_k|\right)^{s-1}\\
  &\le G_n\|\theta\|\sum_{s=1}^\iy \frac{1}{2^{s-1}}=2G_n\|\theta\|.
  \end{align*}
Estimate \eqref{7.14} follows from \eqref{7.13}, the well-known theorem on changing summation order for the series with
positive terms (see \cite[\S1.6]{9}) and \eqref{8.41}:
  \begin{align*}
   \Phi_n(\theta)&=\sum_{k=n}^\iy|f_k(\theta)|=\sum_{k=n}^\iy\frac{1}{r_kv_ku_{k+1}}\left|\sum_{s=1}^\iy(-1)^s(\mathcal A^s
  g(\theta))_k\right|
   \\
    &\le\sum_{k=n}^\iy\frac{1}{r_kv_ku_{k+1}}\sum_{s=1}^\iy|(\mathcal
    A^sg(\theta))_k|=\sum_{s=1}^\iy\sum_{k=n}^\iy\frac{|(\A^sg(\theta))_k|}{r_kv_ku_{k+1}}=\varkappa_n(\theta)\le
2G_n\|\theta\|.
  \end{align*}
\end{proof}

\section{Proof of the Criterion for Solvability of the Narrow Hartman-Wintner Problem}\label{ProofCritNar}

In this section we prove \thmref{thm1.3}.

\renewcommand{\qedsymbol}{}
\begin{proof}[Proof of \thmref{thm1.3}] \ {\it Necessity}.
Suppose that problem \eqref{1.4}--\eqref{1.$*$} is solvable.
This means that condition A) from \secref{double} holds, and to prove the necessity of the hypotheses of \thmref{thm1.3} it
remains to prove that the series $J$ (see \eqref{1.7}) converges (at least, conditionally).
Let us prove this.
Throughout the sequel $\{\beta_n\}_{n=0}^\iy$ denotes a solution of problem \eqref{6.2a}--\eqref{6.3a}.
Recall that by \eqref{5.13}--\eqref{5.14} and \eqref{3.2} we have the following basic equalities:
\begin{equation}\label{5.14a}
\mu_n=-\frac{v_n}{u_n}\ \frac{u_{n+1}}{v_{n+1}}C_{n+1}\beta_{n+1}-\frac{v_n}{u_n}\sum_{k=n+1}^\iy\frac{C_{k+1}
\mu_k}{r_kv_kv_{k+1}},\quad n\ge 0.
\end{equation}
\begin{equation}\label{5.13a}
\mu_n=r_nv_nu_{n+1}\Delta\beta_n,\quad n=0,\qquad \lim_{n\to\iy}\mu_n=0,\quad \lim_{n\to\iy}\beta_n=1.
\end{equation}

{}From \eqref{8.17}--\eqref{8.18} it follows that the system of equalities \eqref{5.14a} and, more precisely, its part for
$n\ge n_0,$ can be written in the form
\begin{equation}\label{8.3a}
\mu=g(\beta)-\A\mu,\qquad \mu=\{\mu_n\}_{n=n_0}^\iy,\qquad \beta=\{\beta_n\}_{n=n_0}^\iy.
\end{equation}

Let us now choose $n_0\gg1$ so that the following inequality holds:
\begin{equation}\label{8.4a}
\|\A\|\le \sup_{n\ge n_0}|C_n|\le\frac{1}{2}
\end{equation}
(see \eqref{7.6}).
Then the operator $\A$ in \eqref{8.3a} is compressing, and therefore we can find $\mu:$
\begin{equation}\label{8.23}
\mu=(E+\mathcal A)^{-1}g(\beta)=g(\beta)+\sum_{k=1}^\iy(-1)^k\mathcal A^kg(\beta).
\end{equation}
Here $E$ is the identity operator.
By \eqref{7.7}, equality \eqref{8.23} can be written coordinatewise (for $n\ge n_0):$
\begin{equation}\label{8.24}
\mu_n=g_n(\beta)+\left(\sum_{k=1}^\iy(-1)^k\mathcal A^kg(\beta)\right)_n=g_n(\beta)+\sum_{k=1}^\iy(-1)^k(\mathcal
A^kg(\beta))_n,\quad n\ge n_0.
\end{equation}
{}From \eqref{8.24}, \eqref{5.14a} and \eqref{5.13a}, we get (for $n\ge n_0)$:
\begin{equation}\label{8.7a}
\Delta\beta_n+\frac{C_{n+1}}{r_nu_nv_{n+1}}\beta_{n+1}=\frac{1}{r_nv_nu_{n+1}}\sum_{s=1}^\iy(-1)^s(\A^sg(\beta))_n.
\end{equation}

Our next goal is to analyze equation \eqref{8.7a} in detail (note that its solution $\{\beta_n\}_{n=n_0}^\iy$ exists as a
solution of problem \eqref{6.2a}--\eqref{6.3a}).
For a given $n\ge n_0,$ let us introduce a sequence $\{Q_n^{(k)}\}_{k=0}^\iy$:
\begin{equation}\label{8.8a}
Q_n^{(k)}=\begin{cases} 1,\  &\text{if}\ k=0\\
\prod\limits_{s=0}^{k-1}(1+\eta_{n+s})^{-1},\   &\text{if}\ k\ge 1,\end{cases}\qquad
\eta_m\doe \frac{C_{m+1}}{r_mu_mv_{m+1}},\quad m\ge n_0.
\end{equation}
The sequence $\{Q_n^{(k)}\}_{k=0}^\iy$ can be defined recurrently:
\begin{equation}\label{8.9a}
Q_n^{(0)}=1,\quad Q_n^{(k+1)}=\frac{Q_n^{(k)}}{1+\eta_{n+k}},\qquad k\ge 0,\ n\ge n_0.
\end{equation}
Note that from \eqref{1.13} and \eqref{2.2} it follows that
\begin{equation}\label{8.10a}
\lim_{m\to\iy}\eta_m=0.
\end{equation}
Hence $n_0$ can be chosen so that $|\eta_n|\le 2^{-1}$ for $n\ge n_0,$ and therefore the sequence $\{Q_n^{(k)}\}_{k=0}^\iy$ is
well-defined.
Let us write the solution $\{\beta_n\}_{n=n_0}^\iy$ of equation \eqref{8.7a} in the form
\begin{equation}\label{8.11a}
\beta_{n+k}=d_kQ_n^{(k)},\qquad n\ge n_0,\ k=0,1,2,\dots
\end{equation}
where the factors $d_k,$\ $k\ge n_0$ are to be determined later.
From \eqref{8.11a} it follows that
\begin{equation}\label{8.12a}
\beta_n=\beta_{n+k}\bigm|_{k=0}=d_0Q_n^{(0)}=d_0\quad\Rightarrow \quad d_0=\beta_n.
\end{equation}
To find $d_k$ for $k>0,$ let us plug \eqref{8.11a} into \eqref{8.7a} with number $(n+k),$\ $n\ge n_0,$\ $k=1,2,\dots$
Here we use the notations \eqref{8.41} and \eqref{8.8a} (see also \eqref{8.9a}):
  \begin{equation}
\begin{aligned}\label{8.13a}
f_{n+k}(\beta)&=\Delta\beta_{n+k}+\frac{C_{n+k+1}}{r_{n+k}u_{n+k}v_{n+k+1}}\beta_{n+k+1}=d_{k+1}Q_n^{(k+1)}-d_kQ_n^{(k)}+\eta_{n+k}
d_{k+1}Q_n^{(k+1)}\\
&=Q_n^{(k+1)}(1+\eta_{n+k})d_{k+1}-d_kQ_n^{(k)}=(d_{k+1}-d_k)Q_n^{(k)}\\
&\Rightarrow\begin{cases}
 d_{k+1}=d_k+\frac{1}{Q_n^{(k)}}f_{n+k}(\beta),\quad k\ge 0,\ n\ge n_0\\
d_0=\beta_n\end{cases}
\end{aligned}
 \end{equation}

{}From \eqref{8.13a} it immediately follows that
\begin{equation}\label{8.14a}
d_{k+1}=\beta_n+\sum_{\ell=0}^k\frac{1}{Q_n^{(\ell)}}f_{n+\ell}(\beta),\qquad k\ge 0, \ n\ge n_0.
\end{equation}
Finally, from \eqref{8.14a} and \eqref{8.11a} we get
\begin{equation}\label{8.15a}
\beta_{n+k+1}=d_{k+1}Q_n^{(k+1)}=\beta_nQ_n^{(k+1)}+\sum_{\ell=0}^k\frac{Q_n^{(k+1)}}{Q_n^{(\ell)}}f_{n+\ell}(\beta),\qquad
k\ge0,\ n\ge n_0.
\end{equation}

Below we study the sequence $\{Q_n^{(k)}\}_{k=0}^\iy $ in detail for $n\ge n_0.$
It turns out that \eqref{8.15a} admits a more convenient form.
To obtain it, we set $k=m-1$,\ $m\ge 1$ in \eqref{8.15a}.
Then
  \begin{equation}\label{8.16a}
   \beta_{n+m}=\beta_nQ_n^{(m)}+\sum_{\ell=0}^{m-1}\frac{Q_n^{(m)}}{Q_n^{(\ell)}}f_{n+\ell}(\beta),\qquad m\ge 1,\ n\ge n_0.
 \end{equation}
Let us now check the equality
  \begin{equation}\label{8.17a}
  Q_n^{(m)}=Q_n^{(\ell)}
 \cdot Q_{n+\ell}^{(m-\ell)},\qquad m\ge 1,\ \ell=0,1,\dots,m-1;\ n\ge n_0.
  \end{equation}

We consider two separate cases:
$$ 1)\quad \ell=0,\quad m\ge 1;\qquad 2)\quad \ell=1,\dots,m-1;\quad m\ge 2.$$

1)\quad If $\ell=0$ and $m\ge 1,$ then \eqref{8.8a} implies
$$\frac{Q_n^{(m)}}{Q_n^{(\ell)}}\bigm|_{\ell=0}=\frac{Q_n^{(m)}}{Q_n^{(0)}}=Q_n^{(m)}=Q_{n+\ell}^{(m-\ell)}\bigm|_{\ell=0}.$$

2)\quad If $\ell=1,\dots,m-1;\ m\ge 2$ then \eqref{8.8a} implies
\begin{align*}
\frac{Q_n^{(m)}}{Q_n^{(\ell)}}&=\frac{(1+\eta_n)\dots (1+\eta_{n+\ell-1})}{(1+\eta_n)\dots(1+\eta_{n+m-1})}=\frac{1}
{(1+\eta_{n+\ell})\dots (1+\eta_{n+m-1})}\\
&=\frac{1}{(1+\eta_{n+\ell})\dots (1+\eta_{n+\ell+(m-\ell-1})}=Q_{n+\ell}^{(m-\ell)}.
\end{align*}

Let $m\ge 1.$ Denote
\begin{equation}\label{8.18a}
M_n^{(m)}=\max_{\ell\in[0,m-1]}\left|\frac{Q_n^{(m)}}{Q_n^{(\ell)}}-1\right|,\quad [0,m-1]\doe\{0,\dots,m-1\}.
\end{equation}
Clearly, there is $\ell_0\in[0,m-1]$ such that (see \eqref{8.17a}):
\begin{equation}\label{8.19a}
M_n^{(m)}=\left|\frac{Q_m^{(m)}}{Q_n^{(\ell_0)}}-1\right|=\left|Q_{n+\ell_0}^{(m-\ell_0)}-1\right|.
\end{equation}
Then from \eqref{8.19a} and \eqref{8.16a} it follows that
\begin{align*}
\beta_{n+m}&=\beta_{n+\ell_0+(m-\ell_0)}=\beta_{n+\ell_0}Q_{n+\ell_0}^{(m-\ell_0)}+\sum_{\ell=0}^{m-\ell_0-1}\frac{Q_{n+\ell_0}^{(
m-\ell_0)}}{Q_{n+\ell_0}^{(\ell)}}f_{n+\ell_0+\ell}(\beta)\\
&=\beta_{n+\ell_0}Q_{n+\ell_0}^{(m-\ell_0)}+\sum_{\ell=0}^{m-\ell_0-1}Q_{n+\ell_0+\ell}^{(m-\ell_0-\ell)}f_{n+\ell_0+\ell}(\beta).
\end{align*}
The last inequality implies
\begin{equation}\label{8.20a}
\begin{aligned}
\beta_{n+\ell_0}[Q_{n+\ell_0}^{(m-\ell_0)}-1]&=\beta_{n+m}-\beta_{n+\ell_0}-\sum_{\ell=0}^{m-\ell_0-1}[Q_{n+\ell_0+\ell}
^{(m-\ell_0-\ell)}
-1]f_{n+\ell_0+\ell}(\beta) \\
& -\sum_{\ell=0}^{m-\ell_0-1}f_{n+\ell_0+\ell}(\beta),\quad m\ge1,\  n\ge  n_0.
\end{aligned}
\end{equation}

Let us replace, if necessary, $n_0$ by a bigger number so that $|\beta_n|\ge 2^{-1}$ for all $n\ge n_0.$
Furthermore, note that \eqref{8.17a} implies
\begin{equation}\label{8.21a}
\begin{aligned}
\max_{\ell\in[0,m-\ell_0-1]}|Q_{n+\ell_0+\ell}^{(m-\ell_0-\ell)}-1|&\le\max_{\ell\in[0,m-1]}|Q_{n+\ell}^{(m-\ell)}-1|=\max_{s\in
[0,m-1]}\left|\frac{Q_n^{(m)}}{Q_n^{(s)}}-1\right|=M_n^{(m)},\\
&\qquad\quad\qquad m\ge 1,\ n\ge n_0.
\end{aligned}
\end{equation}
Therefore from \eqref{8.20a} and \eqref{8.21a}, we get
\begin{align*}
M_n^{(m)}&=|Q_{n+\ell_0}^{(m-\ell_0)}-1|\le 2|\beta_{n+m}-\beta_{n+\ell_0}|+2\sum_{\ell=0}^{m-\ell_0-1}|Q_{n+\ell_0+\ell}^{(m-
\ell_0-\ell)}-1|\ |f_{n+\ell_0+\ell}(\beta)|\\
&+\sum_{\ell=0}^{m-\ell_0-1}|f_{n+\ell_0+\ell}(\beta)|\le 4\max_{k\in[0,m]}|\beta_{n+k}-\beta_n|+2\max_{\ell\in[0,m-\ell_0-1]}
\left|Q_{n+\ell_0+\ell}^{(m-\ell_0-1)}-1\right|\sum_{k=n}^\iy f_k(\beta)|\\
&+2\sum_{k=n}^\iy|f_k(\beta)|\le
4\max_{k\in[0,m]}|\beta_{n+k}-\beta_n|+2M_n^{(m)}\sum_{k=n}^\iy|f_k(\beta)|+2\sum_{k=n}^\iy|f_k(\beta)|.
\end{align*}

Thus for all $m\ge 1$ and $n\ge n_0,$ the following inequality holds:
\begin{equation}\label{8.22a}
M_n^{(m)}\le 4\max_{k\in[0,m]}|\beta_{n+k}-\beta_n|+2M_n^{(m)}\sum_{k=n}^\iy|f_k(\beta)|+2\sum_{k=n}^\iy|f_k(\beta)|.
\end{equation}
From Lemma 7.3 we conclude that the estimate \eqref{8.22a} can be continued as follows (for all $m\ge 1,$\ $n\ge n_0):$
\begin{equation}\label{8.23a}
M_n^{(m)}\le 4\max_{k\in[0,m]}|\beta_{n+k}-\beta_n|+4M_n^{(m)}G_n\|\beta\|+4G_n\|\beta\|.
\end{equation}
Let $\ve$ be a given positive number.
From condition A) (see \S \ref{double}) we conclude that there is $n_0(\ve)\ge n_0$ such that for all $n\ge n_0(\ve)$, the
following inequalities hold:
\begin{equation}\label{8.24a}
\max_{k\in[0,m]}|\beta_{n+k}-\beta_n|\le \frac{\ve}{16},\qquad 4G_n\|\beta\|\le\min\left\{\frac{1}{2},\frac{\ve}{4}\right\}.
\end{equation}
Then from \eqref{8.23a} and \eqref{8.24a} for all $n\ge n_0(\ve),$\ $m\ge 1$, we get
\begin{equation}\label{8.25a}
M_n^{(m)}\le\frac{\ve}{4}+\frac{1}{2}M_n^{(m)}+\frac{\ve}{4}\Rightarrow M_n^{(m)}\le \ve.
\end{equation}
By \eqref{8.8a}, \eqref{8.19a} and \eqref{8.25a}, we now conclude that
$$\left|\frac{Q_n^{(m)}}{Q_n^{(0)}}-1\right|=|Q_n^{(m)}-1|\le M_n^{(m)}
\le\ve\qquad \text{for}\ n\ge n_0(\ve),\ m\ge 1.$$ Now by
\eqref{8.8a},  \eqref{8.19a} and \eqref{8.25a}, we get
$$|Q_n^{(m)}-1|=\left|\frac{Q_n^{(m)}}{Q_n^{(0)}}-1\right|\le M_n^{(n)}\le\ve,
\quad n\ge n_0(\ve), \quad m\ge 1,$$ or, equivalently,
\begin{equation}\label{8.26a}
\left|\prod_{s=0}^m\frac{1}{1+\eta_{n+s}}-1\right|<\ve\qquad\text{for}\ n\ge n_0(\ve),\ m\ge 0.
\end{equation}

\begin{lem}\label{lem8.1}
Let $\{a_n\}_{n=0}^\iy$ be a sequence of (maybe complex) numbers such that $|a_n|\le 2^{-1}$ for all $n\ge 0.$
The infinite product
$$A=\prod_{n=0}^\iy(1+a_n)$$
converges to a finite nonzero number if and only if for every
$\ve>0$ there exists $n_0(\ve)\gg 1$ such that for all $n\ge
n(\ve)$ and all $m\ge 1$, the following inequality holds:
\begin{equation}\label{8.27a}
\left|\prod_{k=1}^m(1+a_{n+k})-1\right|\le \ve.
\end{equation}
\end{lem}

\renewcommand{\qedsymbol}{}
\begin{proof}[Proof of \lemref{lem8.1}] \ {\it Necessity}.
Suppose that the product  $A$ converges to a finite nonzero number.
Denote
\begin{equation}\label{8.28a}
A_n=\prod_{k=1}^n(1+a_k),\qquad
A_n^{(m)}=\prod_{k=1}^m(1+a_{n+k}),\qquad m\ge 1,\ n\ge0.
\end{equation}
Then the sequence $\{A_n\}_{n=0}^\iy$ converges, and by Cauchy's
criterion for every $\ve>0$ there is $n(\ve)$ such that for all
$n\ge n(\ve)$ and $m\ge1$, the following inequality holds:
\begin{equation}\label{8.29a}
|A_{n+m}-A_n|\le\ve\frac{|A|}{2}\quad\Rightarrow\quad |A_n|\ |A_n^{(m)}-1|\le\ve\frac{|A|}{2}.
\end{equation}
Let us replace, if necessary, $n(\ve)$ by a bigger number  so that for all $n\ge n(\ve)$, the following inequality holds:
$$|A_n|\ge \frac{|A|}{2},\qquad n\ge n(\ve).$$
Then from \eqref{8.29a} for $n\ge n(\ve)$ and $m\ge 1$, we get
$$\frac{1}{2}|A|\ |A_n^{(m)}-1|\le |A_n| \ |A_n^{(m)}-1|\le \ve\frac{|A|}{2}\quad\Rightarrow\quad
|A_n^{(m)}-1|\le\ve\quad\text{for}\ n\ge n(\ve),\ m\ge1.$$
\end{proof}

\renewcommand{\qedsymbol}{\openbox}
\begin{proof}[Proof of \lemref{lem8.1}] \ {\it Sufficiency}.
Suppose that \eqref{8.27a} holds.
Then there is a sequence of integers $\{n_k\}_{k=1}^\iy$ such that
\begin{equation}\label{8.30}
|A_n^{(m)}-1|\le \frac{1}{2^k}
\end{equation}
for all $n\ge n_k,$\  $m\ge 1.$ Consider the value $A_n.$ Suppose
that
$$n_k<n\le n_{k+1},\qquad k=1,2,\dots$$
Then
\begin{alignat*}{2}
&A_n &&=A_{n_1}\cdot A_{n_1}^{(n_2-n_1)}\cdot A_{n_2}^{(n_3-n_2)}\dots
A_{n_k}^{(n-n_k)} \quad\Rightarrow\\
&|A_{n_1}|&&=|A_{n_1}|\ |(A_{n_1}^{(n_2-n_1)}-1)+1|\ |(A_{n_2}^{(n_3-n_2)}-1)+1|
\dots|(A_{n_k}^{(n-n_k)}-1)+1|\\
& && \le |A_{n_1}|\left(1+\frac{1}{2
}\right)\left(1+\frac{1}{2^{2}}\right)\dots
\left(1+\frac{1}{2^{k}}\right)\le\tau<\iy,\\
&|A_n|&&=|A_{n_1}|\ |(A_{n_1}^{(n_2-n_1)}-1)+1|\ |(A_{n_2}^{(n_3-n_2)}-1)+1|\dots
 |(A_{n_k}^{(n-n_k)}-1)+1|\\
& &&\ge
|A_{n_1}|\left(1-\frac{1}{2}\right)\left(1-\frac{1}{2^{2}}\right)\dots
\left(1-\frac{1}{2^{k}}\right)\ge \frac{1}{\tau}>0.
\end{alignat*}
Thus
\begin{equation}\label{8.31a}
\tau^{-1}\le|A_n|\le \tau,\qquad n\ge 1.
\end{equation}

Let $\tau$ be the number from inequality \eqref{8.31a}.
Let now $\ve>0$ be given, and choose $n(\ve)$ so that the following inequalities hold:
\begin{equation}\label{8.32a}
|A_n^{(m)}-1|\le\frac{\ve}{\tau}\qquad\text{for}\quad n\ge
n(\ve),\ m\ge1.
\end{equation}
Let us verify that then
\begin{equation}\label{8.33a}
|A_{n+m}-A_n|\le\ve\qquad\text{for all}\ n\ge n(\ve),\ m\ge1.
\end{equation}
We have
$$|A_{n+m}-A_n|=|A_n| \ |A_n^{(m)}-1|\le\tau\frac{\ve}{\tau}=\ve,\quad n\ge n(\ve),
\ m\ge1.$$ Thus the sequence $\{A_n\}_{n=1}^\iy$  converges, and
by \eqref{8.31a} its limit is finite and nonzero, which   was to
be proved.
\end{proof}

Let us  now consider the infinite product
\begin{equation}\label{8.34a}
P=\prod_{n=n_0}^\iy(1+\eta_n),\qquad \eta_n=\frac{C_{n+1}}{r_nu_nv_{n+1}},\qquad n\ge  n_0.
\end{equation}

Recall (see \eqref{8.10a} and below) that $|\eta_n|\le 2^{-1}$ for $n\ge n_0.$
Let us show that from condition A)  it follows that the product $P$ converges to a finite nonzero number.
Suppose (see \eqref{8.8a}) that $P_n^{(m)}$ is defined by the equality
\begin{equation}\label{8.35a}
P_n^{(m)}\cdot Q_n^{(m)}=1,\qquad n\ge n_0,\ m\ge 0.
\end{equation}
According to \eqref{8.26a}, for any $\ve>0$ there is $n(\ve)$ such that for all $n\ge n(\ve)$ and $m\ge0$, the following
inequality holds:
\begin{equation}\label{8.36a}
\left|\frac{1}{P_n^{(m)}}-1\right|<\ve\quad\Rightarrow\quad |P_n^{(m)}-1|<\ve|P_n^{(m)}|.
\end{equation}
Let $\ve=\frac{1}{2}.$
Then there is $n\left(\frac{1}{2}\right)$ such that for all $n\ge n\left(\frac{1}{2}\right)$ and $m\ge 0,$ we have (in view of
\eqref{8.36a}):
$$\left|\frac{1}{P_n^{(m)}}\right|=\left|\frac{1}{P_n^{(m)}}-1+1\right|\ge 1-\left|\frac{1}{P_n^{(m)}}-1\right|\ge
1-\frac{1}{2}=\frac{1}{2}\quad\Rightarrow\quad |P_n^{(m)}|\le 2$$
$$\text{for all}\quad n\ge n\left(\frac{1}{2}\right),\ m\ge 0.$$

Let now $\ve>0$ be an arbitrary number, and choose $m(\ve)\ge n\left(\frac{1}{2}\right)\ge n_0$ so that for all $n\ge n(\ve)$
and $m\ge 0$
\begin{equation}\label{8.37a}
\left|\frac{1}{P_n^{(m)}}-1\right|<\frac{\ve}{2}\quad\Rightarrow\quad
|P_n^{(m)}-1|\le\frac{\ve}{2}|P_n^{(m)}|<\frac{\ve}{2}\cdot2 =\ve.
\end{equation}
{}From \eqref{8.37a} and \lemref{lem8.1}, it follows that the product $P$ (see \eqref{8.34a}) converges to a finite nonzero
number.
This implies that the series $H_{n_0}$ (see \eqref{7.4a})  converges (at least, conditionally).
To check that note that the series
$$B=\sum_{n=0}^\iy\eta_n^2,\qquad \eta_n=\frac{C_{n+1}}{r_nu_nv_{n+1}} ,\qquad n\ge 0$$
absolutely converges.
Indeed, from \eqref{2.2} and condition A) (see \secref{double}), we get
\begin{equation}\label{8.38a}
\sum_{n=n_0}^\iy|\eta_k|^2=\sum_{k=k_0}^\iy\left(\frac{|C_{k+1}|}{r_ku_kv_{k+1}}\right)^2=\sum_{k=n_0}^\iy\frac{1}{
r_ku_kv_{k+1}}\ \frac{|C_{k+1}|^2}{r_ku_kv_{k+1}}\le G<\iy.
\end{equation}

Furthermore, since $\eta_n\to 0$ as $n\to\iy,$ there is $n_0\gg 1$ such that for all $n\ge n_0$ we have
\begin{equation}\label{8.39a}
\Ln(1+\eta_n)=\eta_n+O(|\eta_n|^2),\qquad n\ge n_0.
\end{equation}
where the constant in $O(\cdot)$ is absolute.
Denote
$$S_n=\prod_{k=n_0}^n(1+\eta_k),\qquad B_n=\sum_{k=n_0}^n|\eta_k|^2,\qquad n\ge n_0.$$
Then by \eqref{8.39a} we get
$$\Ln S_n=\sum_{k=n_0}^n\eta_k+O\left(\sum_{k=n_0}^n|\eta_k|^2\right)=\sum_{k=n_0}^n\eta_k+O(B_n),\quad n\ge n_0.$$
Since the sequences $\{\Ln S_n\}_{n=n_0}^\iy,$\ $\{B_n\}_{n=n_0}^\iy$ in the last equality have finite limits, the series
$H_{n_0}$ (see \eqref{7.4a}) converges  (at least, conditionally).
This almost immediately implies that the series $J$ (see \eqref{1.7}) converges.
Indeed, since
\begin{equation}\label{8.40a}
C_{n+1}-C_n=\frac{C_{n+1}}{r_nu_nv_{n+1}}-\sigma_nu_nv_n,\qquad n\ge 0
\end{equation}
(see the proof of \lemref{lem7.1a}), after adding up the equalities \eqref{8.40a}, we obtain
$$C_{m+1}-C_{n_0}=\sum_{n=n_0}^m\frac{C_{n+1}}{r_nu_nv_{n+1}}-\sum_{n=n_0}^m\tau_nu_nv_n,\qquad m\ge n_0.$$
Condition A) (see \S\ref{double}) now implies:
$$-C_{n_0}=\lim_{m\to\infty}(C_{n+1}-C_{n_0})=H_{n_0}-\lim_{m\to\infty}\sum_{n=n_0}^m
\sigma_nu_nv_n,$$ that is, the series $J$ converges (at least,
conditionally), which   was to be proved.

\renewcommand{\qedsymbol}{}
\begin{proof}[Proof of \thmref{thm1.3}] \ {\it Sufficiency}.
Below we need the following simple  auxiliary assertion.
\end{proof}
\begin{lem}\label{lem8.2}
Suppose that condition B) holds (see \S\ref{double}).
Then there is $n_0\gg 1$ such that the infinite product $P$ (see \eqref{8.34a}) converges to a finite nonzero limit.
\end{lem}

 \renewcommand{\qedsymbol}{\openbox}
\begin{proof}
It is known \cite[\S1.43,  Example (1)]{9}, that if the series
\begin{equation}\label{8.41a}
\sum_{n=n_0}^\iy a_n,\qquad \sum_{n=n_0}^\iy |a_n|^2
\end{equation}
converge (the first one at least conditionally), then the infinite product $\prod\limits_{k=n_0}^\iy(1+a_n),$\ $n_0\gg 1$
converges to a finite nonzero limit.
In our case, $a_n=\eta_n,$\ $n\ge0$ (see \eqref{8.8a}).
{}From condition B) it follows that the second series in \eqref{8.41a} converges (see \eqref{8.38a}), and \lemref{lem8.1}
implies that the first series in \eqref{8.41a} converges.
Finally, from the convergence of any series in \eqref{8.41a} we conclude that $|\eta_n|\le 2^{-1}$ for $n_0\gg 1.$
\end{proof}

Let us introdue an operator $T^{(n_0)}$ acting in the space of sequences \eqref{7.1a} by the rule
\begin{equation}\label{8.42a}
T^{(n_0)}(\theta)=\{(T^{(n_0)}(\theta))_n\}_{n=n_0}^\iy,\qquad \|\theta\|<\iy.
\end{equation}
Here
\begin{equation}\label{8.43a}
(T^{(n_0)}(\theta))_n=-\sum_{k=n}^\iy\frac{P_n}{P_k}f_k(\theta),\qquad n\ge n_0.
\end{equation}
\begin{equation}\label{8.44a}
P_n=\prod_{k=n}^\iy(1+\eta_k),\quad n\ge  n_0,\qquad \eta_k=\frac{C_{k+1}}{r_ku_kv_{k+1}},\quad k\ge n
\end{equation}
and the values $\{f_k(\theta)\}_{k=n_0}^\iy$ are defined by relations \eqref{8.41}, \eqref{8.17} and \eqref{8.18}.
Let us verify that under condition B) there exists $n_0\gg1$ such that for all $n\ge n_0$, the following inequalities hold (see
\eqref{8.50}):
\begin{equation}\label{8.45a}
|(T^{(n_0)}\theta)_n|\le\tau G_n\|\theta\|.
\end{equation}
Indeed, if $n_0\gg 1,$ by \lemref{lem8.2} the product $P$ (see \eqref{8.34a}) converges, and therefore $|P_n|\le\tau|P_k|$ for
all $k\ge n\ge n_0.$
Then \eqref{8.45} and \lemref{lem7.3} imply \eqref{8.45a}:
$$|(T^{(n_0)}(\theta)_n|\le\sum_{k=n}^\iy\left|\frac{P_n}{P_k}\right||f_k(\theta)|\le\tau\sum_{k=n}^\iy|f_k(\theta)|=\tau\Phi_n(\theta)
\le\tau G_n\|\theta\|.$$

{}From the proven estimate \eqref{8.45a}, we conclude that $\|T^{(n_0)}\|\le\tau G_{n_0}.$
Choosing a bigger $n_0,$ if necessary,  we get the estimate $\tau G_{n_0}\le 2^{-1}.$
Then we finally obtain
\begin{equation}\label{8.46a}
\|T^{(n_0)}\|\le \frac{1}{2}.
\end{equation}
 \begin{remark}\label{rem8.1}
Recall that according to the proof of \lemref{lem7.3} the validity of inequality \eqref{7.14} is guaranteed by a choice of
$n_0$ such that $\sup\limits_{n\ge n_0}|C_n|\le 2^{-1}.$
Therefore below we assume $n_0$ chosen big enough so that apart from the above requirements to $n_0,$ the following
inequalities hold together:
\begin{equation}\label{8.47a}
A_{n_0}\le\frac{1}{4},\qquad \tau G_{n_0}\le \frac{1}{2}
\end{equation}
(see Assertion \ref{assert1}).
In \eqref{8.47a} $\tau$ stands for the number from estimate \eqref{8.45a}.
\end{remark}

Consider an equation in the space of sequences \eqref{7.1a} which is essential for that which follows:
\begin{equation}\label{8.48a}
\beta=P+T^{(n_0)}\beta.
\end{equation}
Here $n_0$ is chosen so that the estimate \eqref{8.46a} holds, $P=\{P_n\}_{n=n_0}^\iy$ (see \eqref{8.44a}).
Clearly, the operator $T^{(n_0)}$ is compressing, and therefore equation \eqref{8.48a} has a unique solution $\beta$ in the
space \eqref{7.1a}, and
\begin{equation}\label{8.49a}
\|\beta\|\le\tau\|P\|=\tau<\iy.
\end{equation}
Let us study the properties of the solution $\beta=\{\beta_n\}_{n=n_0}^\iy$ of equation \eqref{8.48a}.
{}From \eqref{8.46a} and \eqref{8.48a}, we get
 \begin{equation}\label{8.50a}
 \beta=(E-T^{(n_0)})^{-1}P=P+\sum_{k=1}^\iy(T^{(n_0)})^kP
 \end{equation}
where $E$ is the identity operator.
Writing down \eqref{8.23} coordinatewise, we get (see \lemref{lem7.2}):
\begin{equation}\label{8.51a}
\beta_n=P_n+\sum_{k=1}^\iy[(T^{(n_0)})^kP]_n,\qquad n\ge n_0.
\end{equation}
Let us show that $\beta_n\to1$ as $n\to\iy.$
Denote
$$\omega=\sum_{k=1}^\iy(T^{(n_0)})^{k-1}P=(E-T^{(n_0)})^{-1}P\quad\Rightarrow\quad \|\omega\|\le\tau<\iy.$$
Then  \eqref{8.50a} takes the form
\begin{equation}\label{8.52a}
\beta=P+T^{(n_0)}\omega\quad\Rightarrow\quad \beta_n=P_n+(T^{(n_0)}\omega)_n,\qquad n\ge n_0.
\end{equation}
{}From \eqref{8.45a} it follows that
$$|(T^{(n_0)}\omega)_n|\le\tau G_n\|\omega\|\le\tau G_n ,$$
and therefore
$$\lim_{n\to\iy}\beta_n=\lim_{n\to\iy}(P_n+O(G_n))=\lim_{n\to\iy} P_n=1.$$

Let us now check that the solution $\beta=\{\beta_n\}_{n=n_0}^\iy$ of equation \eqref{8.48a} is also a solution of equation
\eqref{8.7a}.
{}From \eqref{8.48a}  and the definition of the operator
$T^{(n_0)}$, we get
\begin{equation}\label{8.53a}
\beta_n=P_n-\sum_{k=n}^\iy\frac{P_n}{P_k}f_k(\beta),\qquad n\ge n_0.
\end{equation}
Let us now define a sequence $\hat d=\{\hat d_n\}_{n=n_0}^\iy$ by the equality
\begin{equation}\label{8.54a}
\beta_n=\hat d_n\cdot P_n\quad\Rightarrow\quad \hat d_n=1-\sum_{k=n}^\iy\frac{1}{P_k}f_k(\beta),\qquad n\ge n_0.
\end{equation}
{}From \eqref{8.54a} it follows that
\begin{align}
\beta_{n+1}-\beta_n&=\hat d_{n+1}P_{n+1}-\hat d_nP_n=\hat d_{n+1}P_{n+1}-\left(d_{n+1}-\frac{f_n(\beta)}{P_n}\right)P_n
\nonumber\\
&=\hat d_{n+1}P_{n+1}-\hat d_{n+1}P_n+f_n(\beta)=\hat d_{n+1}(P_{n+1}-P_n)+f_n(\beta)\nonumber\\
&=\hat d_{n+1}P_{n+1}(1-1-\eta_n)+f_n(\beta)=-\eta_n\beta_{n+1}+f_n\nonumber\\
&\Rightarrow \Delta
\beta_n+\frac{C_{n+1}}{r_nu_nv_{n+1}}\beta_{k+1}=\frac{1}{r_nv_nu_{n+1}}\sum_{k=1}^\iy(-1)^k(\A^kg(\beta))_n,\quad
n\ge n_0. \label{8.55a}
\end{align}
{}From \eqref{8.55a} we immediately obtain
\begin{equation}\label{8.56a}
r_nv_nu_{n+1}\Delta\beta_n=-\frac{v_n}{u_n}\
\frac{u_{n+1}}{v_{n+1}}C_{n+1}\beta_{n+1}+\sum_{k=1}^\iy(-1)^k(\A^kg(\beta))_n,\qquad n\ge n_0.
\end{equation}

Then from \eqref{8.47a}, Assertion \ref{assert1}, \eqref{8.18}, \eqref{7.5} and \eqref{2.2} we get
$$
r_nv_nu_{n+1}|\Delta\beta_n| \le \frac{u_{n+1}v_n}{u_nv_{n+1}}|C_{n+1}|\ |\beta_{n+1}|+\sum_{k=1}^\iy|(\A^kg(\beta))_n|
\le \tau A_n+\sum_{k=1}^\iy(2A_n)^k\tau\le \tau A_n,$$
$$ A_n\to0\quad\text{as}\quad n\to \iy,$$
or, finally,
\begin{equation}\label{8.57a}
\lim_{n\to\iy}r_nv_nu_{n+1}\Delta\beta_n=0.
\end{equation}

Let us write down the system of equalities \eqref{8.56a} in the vector form
\begin{equation}\label{8.58a}
\mu=g(\beta)+\sum_{k=1}^\iy(-1)^k\A^k(g(\beta))=(E+\A)^{-1}g(\beta).
\end{equation}
In \eqref{8.58a}, $E$ is the identity operator, $\mu=\{\mu_n\}_{n=n_0}^\iy,$\ $\mu_n=r_nu_nv_{n+1}\Delta\beta_n,$\ $n\ge n_0,$\
$g(\beta)$ and the operator $\A$ is defined according to \eqref{8.18} and \eqref{8.17}.
{}From \eqref{8.58a} it follows that
\begin{equation}\label{8.59a}
\mu=g(\beta)-\A\mu.
\end{equation}
Equality \eqref{8.59a} written in the coordinate form looks as follows:
\begin{equation}\label{8.60a}
\mu_n=-\frac{u_{n+1}}{u_n}\ \frac{v_n}{v_{n+1}}C_{n+1}\beta_{n+1}-\frac{v_n}{u_n}\sum_{k=n+1}^\iy\frac{C_{k+1}
\mu_k}{r_kv_kv_{k+1}},\quad n\ge 0.
\end{equation}
Let us plug the values $\mu_n,$\ $n\ge n_0$ and $C_n,$\ $n\ge n_0$  (see \eqref{1.8})  into \eqref{8.60a}.
We obtain \eqref{6.2a}:
$$r_nu_nu_{n+1}\Delta
 \beta_n=-\beta_{n+1}\sum_{k=n+1}^\iy\sigma_ku_k^2-\sum_{k=n+1}^\iy\Delta\beta_k\sum_{m=k+1}^\iy\sigma_mu_m^2,\qquad n\ge n_0.$$
 Thus, under the hypotheses of \thmref{thm1.3} there exists a solution of problem \eqref{6.2a}--\eqref{6.3a}, and therefore
problem \eqref{1.4}--\eqref{1.6} is solvable (see \lemref{lem6.2}).
Since the series $G$ converges by assumption, from \thmref{thm1.4}  it  follows that the narrow Hartman-Wintner problem is also
solvable.
\end{proof}

\begin{proof}[Proof of \corref{cor1.3.1}]

I)\ Suppose that the series $J$ (see \eqref{1.7}) converges (at least, conditionally).
{}From \eqref{1.15} and \eqref{1.9} it follows that
$$
\sum_{n=0}^\iy\frac{|\Re(J_{n+1}\ov C_{n+1})|}{r_nu_nv_{n+1}}\le\sum_{n=0}^\iy\frac{|J_{n+1}|\
|C_{n+1}|}{r_nu_nv_{n+1}}\le \sum_{n=0}^\iy\frac{A_{n+1}|C_{n+1}|}{r_nu_nv_{n+1}}.$$
Therefore, if any of the first two inequalities of \eqref{1.16} holds, the series $L$ (see \eqref{1.15}) absolutely
converges, and by  \thmref{thm1.3}  problem
\eqref{1.4}--\eqref{1.$*$} is solvable.
Furthermore,  if the third series in \eqref{1.16} converges, by \lemref{lem5.2} the series
$G$ also converges, and problem
\eqref{1.4}--\eqref{1.$*$} is solvable in view of \thmref{thm1.3}.
\end{proof}

\renewcommand{\qedsymbol}{\openbox}
\begin{proof}[Proof of \thmref{thm1.1}]
This statement immediately follows from the proven part I) of \corref{corl3.1.1}.
Here assumption \eqref{1.11} is not used and is superfluous.

II)\ If the series $J$  (see \eqref{1.7}) and $P$ (see \eqref{1.17}) converge (at least, conditionally), by \lemref{lem1.1} the
seres $G$ (see \eqref{1.14}) also converges, and then by \thmref{thm1.3} problem \eqref{1.4}--\eqref{1.$*$} is solvable.

III)\ If the series $J$ (see \eqref{1.7}) and the series \eqref{1.18} converge, then the series $P$ absolutely converges, and
problem \eqref{1.4}--\eqref{1.$*$} is solvable by part II) of this corollary.

IV)\ If the series $J$ absolutely converges, inequality \eqref{1.18} holds, and problem \eqref{1.4}--\eqref{1.$*$} is solvable
by part III) of this corollary.
  \end{proof}

\section{Examples}\label{examples}
In this section we present three examples illustrating the power of Theorems \ref{thm1.2}, \ref{thm1.3} and \ref{thm1.4} and
the simplest tools for their application.
For the reader's convenience and for brevity we analyze all the examples according to a common scheme consisting of steps  I-IV.

\noindent{I.\ Study of the FSS $\{u_n,v_n\}_{n=n_0}^\iy$ of the basic equation \eqref{1.2}}.

In the examples below, the FSS  of equation \eqref{1.2} can be easily found in an explicit form.
To apply Theorems \ref{thm1.2}, \ref{thm1.3} and \ref{thm1.4} in these examples,  it is enough to have sharp by order estimates
of the FSS.
Here such estimates can be obtained by standard methods without any difficulties, using the classic Cauchy-Maclaurain Theorem
\cite[vol. II,  \S2.373]{5}, and therefore we do not prove them.

\noindent{II.\ Study of necessary conditions for the solvability of problem \eqref{1.4}--\eqref{1.6}}.

In this step such conditions are obtained with the help of \thmref{thm1.2}.
More precisely, we find conditions for the convergence of the series $\sigma$ (see \eqref{1.12}), present a priori sharp by
order estimates of the terms  of the sequence $\{C_n\}_{n=n_0}^\iy$ (see \eqref{1.8}), and, as a corollary, precise conditions
guaranteeing equality \eqref{1.13}.

\noindent{III.\ Study of the condition for the coincidence of problems \eqref{1.4}--\eqref{1.6} and \eqref{1.4}--\eqref{1.$*$}}.

In this step, following \thmref{thm1.4}, we establish precise conditions for the convergence of the series $G$ (see
\eqref{1.14}).

\noindent{IV.\ Study of conditions for the solvability of problem \eqref{1.4}--\eqref{1.$*$}}.

In  this step, following \thmref{thm1.3}, we find precise conditions for the convergence of the series $J$ (see \eqref{1.7});
conditions for the  convergence of the series $G$ are given in step III.

\medskip

Note that taking into account the above description, below we give  no further comments on the problems arising at each step.
The results obtained in each step  of the above program and the corresponding conclusions are given in the course of the
exposition.
Finally, throughout the sequel the notation $\vp_n\asymp \psi_n,$\ $n\ge n_0$ means that the following inequalities hold:
$$\tau^{-1}|\vp_n|\le|\psi_n|\le\tau|\vp_n|,\qquad n\ge n_0.$$

In the following example, we describe the situation where problem \eqref{1.4}--\eqref{1.6} and \eqref{1.4}--\eqref{1.$*$} do
not coincide, and therefore the analysis of the Hartman-Wintner problem based  on our theorems contains inevitable gaps.

\begin{examp}\label{examp9.1}
Consider problem \eqref{1.4}--\eqref{1.6} for the equations
\begin{equation}\label{9.1}
\Delta(e^{-n}\Delta y_n)=(n+1)^\gamma e^{-(n+1)}y_{n+1},\qquad n\ge 0
\end{equation}
\begin{equation}\label{9.2a}
\Delta(e^{-n}\Delta z_n)=0\cdot z_{n+1},\qquad n\ge 0.
\end{equation}
\end{examp}

In connection with \exampref{examp9.1}, we have the following result.

\begin{assertion}\label{assert9.1}
For equations \eqref{9.1}--\eqref{9.2a}, problem \eqref{1.4}--\eqref{1.6}

1) is not solvable for $\gamma\ge0$ and $\gamma\in \left[-1,-\left.\frac{1}{2}\right)\right.;$

2) is solvable for $\gamma\in(-\iy,-1).$
\end{assertion}

\begin{remark}\label{rem9.1}
For $\gamma\in\left.\left[-\frac{1}{2}\right. ,0\right),$ problem \eqref{1.4}--\eqref{1.6}  is not equivalent to problem
\eqref{1.4}--\eqref{1.$*$}, and therefore the problem on solvability of problem \eqref{1.4}--\eqref{1.6} remains open in this
case.
\end{remark}

\begin{proof}[Proof of Assertion \ref{assert9.1}]
\
I.\ The sequence $u_n\equiv 1,$ \ $n\ge0,$ is a principal solution of \eqref{9.2a} (see \eqref{1.3}).
By formula \eqref{3.7}, we find a non-principal solution $\{v_n\}_{n=1}^\iy$ of equation \eqref{9.2a}:
$$v_{n+1}=u_{n+1}\sum_{k=0}^n\frac{1}{r_ku_ku_{k+1}}=\sum_{k=0}^n e^k\asymp e^{n+1},\qquad n\ge 0.$$

II.\ The series $\sigma$ is of the form (see \eqref{1.12}):
$$\sigma=\sum_{n=1}^\iy\sigma_nu_n^2=\sum_{n=1}^\iy n^\gamma e^{-n}$$
and therefore converges for every $\gamma\in R.$
Consider $C_n,$\ $n\ge 1$ (see \eqref{1.8})
\begin{equation}\label{9.3a}
C_n=\frac{v_n}{u_n}\sum_{k=n}^\iy \sigma_nu_n^2\asymp e^n\sum_{k=n}^\iy k^\gamma e^{-k}\asymp n^\gamma
\end{equation}
Hence $C_n\to0$ as $n\to\iy$ if and only if $\gamma<0.$
This implies that for $\gamma\ge0$   problem \eqref{1.4}--\eqref{1.6}  is not solvable.

III.\ The study of the series $G$ (see \eqref{1.14}) is based on the application of \eqref{9.3a}:
$$G=\sum_{n=1}^\iy\frac{|C_{n+1}|^2}{r_nu_nv_{n+1}}\asymp \sum_{n=1}^\iy n^{2\g}.$$
Hence the series $G$ converge if and only if $2\gamma<-1.$

We conclude that for $\g\in\left[-\left.\frac{1}{2}\right. ,0\right)$ problems \eqref{1.4}--\eqref{1.6} and
\eqref{1.4}--\eqref{1.$*$} are not equivalent, and for $\g\in\left(-\iy,-\left.\frac{1}{2}\right)\right.$ these problems
coincide.

IV.\ By the results of Step I, we have (see \eqref{1.7}) $J=\sum_{n=1}^\iy\sigma_nu_nv_n\asymp \sum_{n=1}^\iy n^\g.$

We conclude that for $\g\in\left[\left.-1,
\frac{1}{2}\right)\right.$ problem \eqref{1.4}--\eqref{1.6} is not
solvable (since the equivalent problem \eqref{1.4}--\eqref{1.$*$}
is not solvable ), and for $\g<-1$ it is solvable (since the
equivalent problem \eqref{1.4}--\eqref{1.$*$} is solvable.
\end{proof}

In  the  next example we consider a situation where problems \eqref{1.4}--\eqref{1.6} and \eqref{1.4}--\eqref{1.$*$} coincide
provided at least one of them is solvable.
Thus it turns out that problem \eqref{1.4}--\eqref{1.6} has been  fully investigated because of the study of problem
\eqref{1.4}--\eqref{1.$*$}.

\begin{examp}\label{examp9.2}
Consider the problem \eqref{1.4}-\eqref{1.6} for the equations
\begin{equation}\label{9.4a}
\Delta(n^\alpha\Delta y_n)=\frac{1}{(n+1)^\beta}y_{n+1},\qquad n\ge 1
\end{equation}
\begin{equation}\label{9.5a}
\Delta(n^\alpha\Delta z_n)=0\cdot z_{n+1},\qquad n\ge1.
\end{equation}
\end{examp}
Here $\alpha\ge0,\ \beta\in R$.

In connection with \exampref{examp9.2}, we have the following
result.

\begin{assertion}\label{assert9.2}
Problem \eqref{1.4}--\eqref{1.6} for equations \eqref{9.4a}--\eqref{9.5a} is solvable if and only if
\begin{equation}\label{9.6a}
\alpha+\beta>2.
\end{equation}
\end{assertion}

\begin{proof}[Proof of Assertion \ref{assert9.2}]
Throughout the sequel, each of the steps I--IV is subdivided into
three substeps according to the value of $\alpha:$\
$\alpha\in[0,1),$\ $\alpha=1,$\ $\alpha\in(1,\iy)$. Such a
subdivision is necesary because the form of the FSS
$\{u_n,v_n\}_{n=1}^\iy$ of equation \eqref{9.5a} changes with
changing the parameter $\alpha\in[0,\iy).$

\smallskip
\noindent I.1)\ $\alpha\in[0,1).$

The sequence $u_n\equiv 1,$\ $n\ge0$ is the principal solution of \eqref{9.5a}.
The non-principal solution $v_n,$ \ $n=1,2,\dots$  is of the form
$$v_{n+1}=u_{n+1}\sum_{k=1}^n \frac{1}{r_ku_ku_{k+1}}=\sum_{k=1}^n\frac{1}{k^\alpha}\asymp (n+1)^{1-\alpha}.$$

\noindent I.2)\ $\alpha=1.$

The sequence $u_n\equiv 1,$\ $n\ge0$ is the principal solution of \eqref{9.5a}.
The non-principal solution $v_n,$ \ $n=1,2,\dots$ is of the form
$$v_{n+1}=u_{n+1}\sum_{k=1}^n\frac{1}{r_ku_ku_{k+1}}=\sum_{k=1}^n \frac{1}{k}\asymp \Ln(n+1).$$

\noindent I.3)\ $\alpha>1.$

The sequence $v_n\equiv 1,$\ $n\ge0$ is the non-principal solution of \eqref{9.5a}, the principal solution $u_n,$ \ $n\ge1$ is
of the form (see \eqref{2.1}):
$$u_n=v_n\sum_{k=n}^\iy\frac{1}{r_kv_kv_{k+1}}=\sum_{k=n}^\iy\frac{1}{k^\alpha}\asymp \frac{1}{n^{\alpha-1}}.$$

\noindent II.1)\ $\alpha\in [0,1).$

In this case the series $\sigma$ is of the form
$$\sigma=\sum_{n=1}^\iy\sigma_nu_n^2=\sum_{n=1}^\iy \frac{1}{n^\beta}<\iy\quad\Leftrightarrow\quad \beta>1.$$

Consider the sequence $C_n,$ $n\ge 1.$
We get
$$C_n=\frac{v_n}{u_n}\sum_{k=n}^\iy\sigma_nu_k^2\asymp  n^{1-\alpha}
\sum_{k=n}^\iy\frac{1}{k^\beta}\asymp
\frac{n^{1-\alpha}}{n^{\beta-1}}=n^{2-\alpha-\beta}.$$

Therefore $C_n\to0$ as $n\to\iy\quad\Leftrightarrow\quad \alpha+\beta>2.$
The collection of the given and the obtained necessary conditions $\alpha\in [0,1),$\ $\beta>1,$\ $\alpha+\beta>2$ is
equivalent to the collection of conditions $\alpha\in[0,1)$,\ $\alpha+\beta>2.$

\smallskip
\noindent II.2)\ $\alpha=1.$

In this case the series $\sigma$ is of the form
$$\sigma=\sum_{k=1}^\iy \sigma_ku_k^2=\sum_{k=1}^\iy \frac{1}{k^\beta}<\iy\quad\Leftrightarrow\quad \beta>1.$$

Consider the sequence $C_n,$\ $n\ge 1.$
We get
$$C_n=\frac{v_n}{u_n}\sum_{k=n}^\iy\sigma_k u_k^2\asymp \sum_{k=n}^\iy\frac{1}{k^\beta}\asymp \frac{\Ln n}{n^{\beta-1}}\to
0\qquad \text{as}\quad n\to 0,$$
since $\beta>1.$
Thus the collection of the given and obtained neccessary conditions $\alpha=1,$\ $\beta>1$ is equivalent to the collection of
conditions $\alpha=1$, \ $\alpha+\beta>2.$

\smallskip
\noindent II.3)\ $\alpha>1.$

In this case the series $\sigma$ is of the form
$$\sigma=\sum_{n=1}^\iy \sigma_nu_n^2\asymp \sum_{n=1}^\iy\frac{1}{n^\beta}\
\frac{1}{n^{2\alpha-2}}<\iy\quad\Leftrightarrow\quad
\beta>3-2\alpha.$$

Consider the sequence $C_n,$\ $n\ge 1.$
We get
$$C_n=\frac{v_n}{u_n}\sum_{k=n}^\iy \sigma_ku_k^2\asymp n^{\alpha-1}\sum_{k=n}^\iy\frac{1}{k^\beta}\
\frac{1}{k^{2\alpha-2}}\asymp
\frac{n^{\alpha-1}}{n^{\beta+2\alpha-3}}\asymp n^{2-\alpha-\beta}\to 0$$
as $n\to\iy\quad\Leftrightarrow\quad\alpha+\beta>2.$
Since for $\alpha>1$ we have $\beta>2-\alpha>3-2\alpha,$ we finally get the collection of necessary conditions in this case:
$\alpha>1,$ \ $\alpha+\beta>2.$

We conclude that if problem \eqref{1.4}--\eqref{1.6} is solvable, then condition \eqref{9.6a} must hold regardless of the value
of $\alpha\in [0,\iy).$

\smallskip
\noindent III.1)\ $\alpha\in[0,1).$

In this case we get the series $G$ in the form
$$G=\sum_{n=1}^\iy \frac{|C_{n+1}|^2}{r_nu_nv_{n+1}}\asymp  \sum_{n=1}^\iy \frac{n^{4-2\alpha-2\beta}}{n^\alpha\cdot 1\cdot
n^{1-\alpha}}=\sum_{n=1}^\iy \frac{1}{n^{2\alpha+2\beta-3}}<\iy\quad\Leftrightarrow\quad \alpha+\beta>2.$$

\noindent III.2)\ $\alpha=1.$

In this case we get the series $G$ in the form
$$G=\sum_{n=1}^\iy\frac{|C_{n+1}|^2}{r_nu_nv_{n+1}}
\asymp \sum_{n=1}^\iy\left(\frac{\Ln n}{n^{\beta-1}}\right)^2\frac{1}{n\Ln n}=\sum_{n=1}^\iy \frac{\Ln
n}{n^{2\beta-1}}<\iy\quad\Leftrightarrow\quad \beta>1.$$
Here the conditions $\alpha=1,$\ $\beta>1$ are equivalent to the conditions $\alpha=1,$\ $\alpha+\beta>2.$

\smallskip
\noindent III.3)\ $\alpha>1.$

In this case we get the series $G $ in the form
$$G=\sum_{n=1}^\iy\frac{|C_{n+1}|^2}{r_nu_nv_{n+1}}\asymp \sum_{n=1}^\iy\frac{1}{n^{2\beta+2\alpha-4}}\ \frac{1}{n^\beta}\cdot
\frac{1}{n^{\alpha-1}}=\sum_{n=1}^\iy\frac{1}{n^{3\alpha+3\beta-5}}<\iy\quad\Leftrightarrow\quad\alpha+\beta>2.$$

We conclude that problem \eqref{1.4}--\eqref{1.6} is solvable  if  and only if problem \eqref{1.4}--\eqref{1.$*$} is solvable.
Here condition \eqref{9.6a} is a necessary condition for solvability of problem \eqref{1.4}--\eqref{1.$*$}.

\smallskip
\noindent IV.1)\ $\alpha\in[0,1).$

In this case the series $J$ is of the form
$$J=\sum_{n=1}^\iy \sigma_nu_nv_n\asymp  \sum_{n=1}^\iy
\frac{n^{1-\alpha}}{n^\beta}=\sum_{n=1}^\iy \frac{1}{n^{\alpha+\beta-1}}<\iy\quad\Leftrightarrow\quad \alpha+\beta >2.$$

\noindent IV.2)\ $\alpha=1.$

In this case the series $J$ is of the form
$$J=\sum_{n=1}^\iy\sigma_nu_nv_n\asymp \sum_{n=1}^\iy\frac{\Ln n}{n^\beta}<\iy\quad\Leftrightarrow\quad
\beta>1\quad\Leftrightarrow\quad\alpha+\beta>2.$$

\newpage

\noindent IV.3)\ $\alpha>1.$

In this case the series $J$ is of the form
$$J=\sum_{n=1}^\iy\sigma_n u_nv_n\asymp \sum_{n=1}^\iy\frac{1}{n^\beta}\cdot\frac{1}{n^{\alpha-1}}=\sum_{n=1}^\iy
\frac{1}{n^{\alpha+\beta-1}}<\iy\quad\Leftrightarrow\quad \alpha+\beta>2.$$

We conclude that the narrow Hartman-Wintner problem is solvable if
and only if condition \eqref{9.6a} holds.
\end{proof}

In the next example the perturbation $\{\sigma_n\}_{n=1}^\iy$ is oscillating, and its absolute value $\{\sigma_n\}_{n=1}^\iy$
coincides with the perturbation from \exampref{examp9.2}.
Thus we show that the exact account of the oscillation of the perturbation allows  one to significantly weaken requirements of
solvability of problem \eqref{1.4}--\eqref{1.6}.

\begin{examp}\label{examp9.3}
Consider the problem \eqref{1.4}--\eqref{1.6} for the equations
\begin{equation}\label{9.7a}
\Delta(n^\alpha\Delta y_n)=\frac{(-1) ^n}{(n+1)^\beta}y_{n+1},\qquad n\ge 1
\end{equation}
\begin{equation}\label{9.8a}
\Delta(n^\alpha\Delta z_n)=0\cdot z_{n+1},\qquad n\ge 1.
\end{equation}
\end{examp}
Here $\alpha\in[0,1)$,\ $\beta\in R.$

In connection with \exampref{examp9.3} we have the following result.

\begin{assertion}\label{assert9.3}
The problem \eqref{1.4}--\eqref{1.6} for equations \eqref{9.7a}--\eqref{9.8a} is solvable if and only if
\begin{equation}\label{9.9a}
\alpha+\beta>1,\quad \beta>0.
\end{equation}
\end{assertion}

\begin{proof}[Proof of Assertion \ref{assert9.3}]
Below in all the steps I--IV, we assume that $\alpha\in[0,1)$ even when it is not specially mentioned.

\noindent I.\ The sequence $u_n\equiv 1,$\ $n\ge 0,$ is the principal solution of \eqref{9.8a}.
The non-principal solution $v_n,$ \ $n\ge 1,$ is of the form
\begin{equation}\label{9.10a}
v_{n+1}=u_{n+1}\sum_{k=1}^n\frac{1}{r_ku_ku_{k+1}}=\sum_{k=1}^n\frac{1}{k^\alpha}\asymp (n+1)^{1-\alpha}.
\end{equation}

\noindent II.\ In this case the series $\sigma$ is of the form
$$\sigma=\sum_{n=1}^\iy\sigma_n u_n^2=\sum_{n=1}^\iy \frac{(-1)^n}{(k+1)^\beta}.$$
{}From the necessary condition for convergence of the series and Leibnitz's theorem, it follows that the series $\sigma$
converges if and only if $\beta>0.$

Consider the sequence $C_n,$ \ $n\ge 1.$
We get
$$C_n=\frac{v_n}{u_n}\sum_{k=n}^\iy \sigma_k u_k^2\asymp n^{1-\alpha}R_n(\beta),\qquad
R_n(\beta)\doe\sum_{k=n}^\iy\frac{(-1)^k}{(k+1)^\beta}.$$
Below we need the following   lemma.

\begin{lem}\label{lem9.1}
For any $\beta>0$ there exists $n_0(\beta)$ such that for all $n\ge n_0(\beta)$, we have
\begin{equation}\label{9.11a}
R_n(\beta)\asymp n^{-\beta}.
\end{equation}
\end{lem}

\begin{proof} Let  $m$ be  any even number.
Then
\begin{align*}
&(-1)^nR_n(\beta) =
 \frac{1}{(n+1)^\beta}-\frac{1}{(n+2)^\beta}+\dots+\frac{1}{(n+m-1)^\beta}-
 \frac{1}{(n+m)^\beta}+(-1)^n
 R_{n+m}(\beta)
 \\
 &=\frac{1}{(n+1)^\beta}\left[1-\left(1-\frac{1}{n+2}\right)^\beta\right]+\dots+\frac{1}{(n+m-1)^\beta}\left[1-\left(1-\frac{1}
 {n+m}\right)
 ^\beta\right]+(-1)^n R_{n+m}(\beta).
\end{align*}
By Taylor's formula, for any $\beta>0$ there exists $n_0(\beta)$ such that for all $k\ge n_0(\beta)$, we have
$$\left(1-\frac{1}{k}\right)^\beta=1-\frac{\beta}{k}+O\left(\frac{1}{k^2}\right),\qquad k\ge n_0(\beta).$$
Here the constant $O(\cdot)$ is absolute and depends only on $\beta.$
Let us now continue the calculation:
$$(-1)^nR_n(\beta)=\beta\sum_{k=1}^{m-1}\frac{1}{(n+k)^\beta(n+k+1)}+O\left(\sum_{k=1}^{m-1}\frac{1}{(n+k)^\beta}\frac{1}{(n+k+1)
^2}\right)+(-1)^nR_{n+m}(\beta).$$
All three summands in the last equality have finite limits as $m\to\iy.$ Hence
$$(-1)^nR_n(\beta)=\beta\sum_{k=1}^\iy\frac{1}{(n+k)^\beta(n+k+1)}+O\left(\sum_{k=1}^\iy\frac{1}{(n+k)^\beta(n+k+1)^2}\right).$$
The obtained equality  and Cauchy-Maclaurin's theorem imply (see \cite[vol.2, \S2.373]{5}) imply \eqref{9.11a}
\end{proof}

Thus for $n\ge n_0(\beta),$ we obtain for the value of $C_n:$
$$C_n\asymp n^{1-\alpha}R_n(\beta)\asymp n^{1-\alpha-\beta}\to0\qquad \text{as}\quad n\to\iy\quad\Leftrightarrow\quad
\alpha+\beta>1.$$ Thus if problem \eqref{1.4}--\eqref{1.6} is solvable, \eqref{9.9a} holds.

\noindent III.\ In this case, for a given $\beta>0,$ for $n\ge n_0(\beta)$ (see \lemref{lem9.1}) we get (see \eqref{8.50}):
\begin{align*}
G_{n_0}(\beta)&=\sum_{n=n_0(\beta)}^\iy\frac{|C_{n+1}|^2}{r_nu_nv_{n+1}}\asymp\sum_{n=n_0(\beta)}^\iy
\frac{1}{n^{2\alpha+2\beta-2}}\ \frac{1}{n^\alpha\cdot n^{1-\alpha}}\\
&=\sum_{n=n_0(\beta)}^\iy \frac{1}{n^{2\alpha+2\beta-1}}<\iy\quad\Leftrightarrow\quad\alpha+\beta>1.
\end{align*}
We conclude that problem \eqref{1.4}--\eqref{1.6} is solvable if and only if problem \eqref{1.4}--\eqref{1.$*$} is solvable,
and condition \eqref{9.9a} is necessary for solvability of problem \eqref{1.4}--\eqref{1.6}.

\noindent IV.\ In this case, for a given $\beta>0$ and the number $n_0(\beta)$ (see \lemref{lem9.1}) we get (see \eqref{1.9},
\eqref{9.10a}):
$$J_{n_0}(\beta)=\sum_{n=n_0(\beta)}^\iy\sigma_n
u_nv_n=\sum_{n=n_0(\beta)}^\iy\frac{(-1)^n}{(n+1)^\beta}\sum_{k=1}^{n-1}\frac{1}{k^\alpha}.$$
Note that under condition \eqref{9.9a} we have
$$\lim_{n\to\iy}R_n(\beta)v_n=0.$$
Indeed, from \eqref{9.10a} and \eqref{9.11a} it follows (for $n\ge n_0(\beta))$ that
$$R_n(\beta)v_n\asymp n^{1-\alpha-\beta}\to0\qquad \text{as}\quad n\to\iy\qquad\text{and}\qquad \alpha+\beta>1.$$
Let us now return to $J_{n_0(\beta)}:$
\begin{align}
J_{n_0(\beta)}&=\sum_{n=n_0(\beta)}^\iy (R_n(\beta)-R_{n+1}(\beta))v_n=\sum_{n=n_0}^\iy[R_n(\beta)v_n-R_{n+1}(\beta)v_{n+1}+R_{n
+1}(\beta)(v_{n+1}-v_n)]\nonumber\\
&=R_{n_0(\beta)}v_{n_0(\beta)}+\sum_{n=n_0(\beta)}^\iy\frac{R_{n+1}(\beta)}{n^\alpha}.\label{9.12a}
\end{align}
But from \eqref{9.9a} and \eqref{9.11a} it follows that the series in the right-hand side of \eqref{9.12a}  absolutely
converges:
$$\sum_{n=n_0(\beta)}^\iy\frac{|R_{n+1}(\beta)|}{n^\alpha}\asymp\sum_{n=n_0(\beta)}^\iy\frac{1}{n^{\alpha+\beta}}<\iy.$$
We conclude that problem \eqref{1.4}--\eqref{1.$*$}  is solvable if \eqref{9.9a} holds.
\end{proof}

We can now prove the assertion stated in \secref{introduction} (see \exampref{examp1}).
Let $\alpha\in [0,1)$ and $1-\alpha<\beta\le\frac{3}{2}-\alpha.$
Then $\beta>0,$\ $\alpha+\beta>1$ and problem \eqref{1.4}--\eqref{1.6} for equations \eqref{9.7a} and \eqref{9.8a} is solvable.
Let us verify that if,  in addition, $\alpha+\beta\le \frac{3}{2},$ then the series B (see \eqref{1.17a}) diverges.
Indeed, from \eqref{9.10a} it follows that
$$B=\sum_{n=1}^\iy [|\sigma_n|u_nv_n]^2\asymp \sum_{n=1}^\iy\frac{n^{2-2\alpha}}{n^{2\beta}}=\sum_{n=1}^\iy
\frac{1}{n^{2\beta+2\alpha-2}}=\iy$$
for $2\alpha+2\beta-2\le 1\quad\Rightarrow\quad \alpha+\beta\le  \frac{3}{2},$ which was to be proved.
\hfill \qed

\end{document}